\documentclass[12pt, titlepage]{article}
\usepackage[sumlimits]{amsmath,mathtools}

\usepackage{amssymb,amscd,amsthm}
\usepackage{enumitem} 
\usepackage[dvipsnames]{xcolor}
\usepackage{graphicx}

\usepackage[utf8]{inputenc} 

\usepackage{setspace}

\usepackage[dvipsnames]{xcolor} 

\newtheorem{theo}{Theorem}
\newtheorem{lem}[theo]{Lemma}

\newtheorem{cor}[theo]{Corollary}

\newtheorem{prop}[theo]{Proposition}

\def\C{{\mathbb C}}

\def\N{{\mathbb N}}
\def\R{{\mathbb R}}
\def\Z{{\mathbb Z}}
\renewcommand{\Pr}{\mathbb P}
\def\eps{\varepsilon}
\def\ph{\varphi}

\def\Om{{\Omega}}
\def\th{{\theta}}
\def\T{{\cal T}}
\def\1{{\mathbf 1}}
\def\d{{\ \rm d}}
\def\osc{\operatorname{osc}}

\newcommand{\norm}[1]{\left\| #1 \right\|}
\newcommand{\abs}[1]{\left| #1 \right|}
\def\B{\mathbf B}
\def\lin{\mathcal L}
\def\bor{\mathcal B}

\DeclareMathOperator*{\supess}{ess\,sup}
\DeclareMathOperator*{\infess}{ess\,inf}
\def\supp{\operatorname{supp}}

\def\off#1{ {\ } \vskip 2 cm}

\def\bbb{{\cal B}}

\begin{document}



\begin{center}{\Large
Random additive
perturbation of a $k$ term 
recurrence relation}\\
\smallskip
{Lisette {\scshape Jager}\footnote[1]{Laboratoire de mathématiques, FR CNRS 3399, EA 4535, Université de Reims Champagne-Ardenne, Moulin de la Housse, BP 1039, 51687 Reims, France}, Killian {\scshape Verdure}}\footnotemark[1] \\
\smallskip
\today
\end{center}

keywords 
: averaged Frobenius Perron operators, dynamical systems, Lasota Yorke inequality, spectral gap, invariant measure

\bigskip

\tableofcontents

\section{Introduction}
 
We aim at studying a class of bounded, real valued  processes $X=(X_t)_{t \ge 0}$
  given by a $k$-term recurrence relation 
\begin{equation}\label{pb1}
  X_{n+k}= \Phi_0 (X_n, \dots, X_{n-k+1}) + \Theta_n,
\end{equation}
where $\Phi_0 $ is a real,  nonlinear, deterministic function defined on
$\R^k$, $X_n : (\Xi, \T, \Pr) \to \R$ are real-valued random variables and 
$ \Theta_n : (\Xi, \T, \Pr) \to \R$ is a stochastic perturbation. We assume
that $\Theta = (\Theta_n)_{n\in \N}$ is i.i.d. and that the distribution
$\Pr_{\Theta_0} $ of one of the variables $\Theta_n$  is absolutely
continuous with respect to the Lebesgue measure $ m$, with density $H$.
\\
When $\Theta$ is equal to $0$, the recurrence relation is deterministic.
Such problems are the subject of \cite{JMN1,JMN2}. In these articles, the
recurrence relation
is imbedded in $\R^k$ in such a way as to give rise to a piecewise
expanding transformation of $\R^k$. The properties of the transformation and
of its transfer operator imply results about the original process $X$.

Expanding maps and their ergodic properties have a long story.
In dimension $1$, \cite{HK, LY73}  studied the transfer
operator of piecewise expanding transformations of an interval, applying
the result of \cite{ITM} to  prove its quasicompacity and the existence
of a spectral gap. This proved very fruitful, but in higher dimensions,
the method does not apply directly.
For example, 
\cite{SOC} studies piecewise expanding transformations  on a hypercube
which is partitioned in smaller hypercubes. In  \cite{TSU1} 
the subdomains are polyhedral and the maps, linear. The domains are not so
regular in \cite{Liv, SAU}, but the (nonlinear)
transformations satisfy very strict and complicated  conditions.
In \cite{YOU}, the
transformation and its properties are presented in a very  axiomatized way, 
which sheds a new light on the transformations and on the subdomains.
But it is extremely difficult to apply to a concrete example,
although \cite{AFLV} proved that such a decomposition in ``towers''
is equivalent to good rates of decay of the correlations. Hence, in situations
where the rate of decay is known, the towers must exist.\\
The difficulty one faces in dimension superior to $1$ is intrinsic: \cite{TSU2} constructs
an example of a piecewise expanding map in dimension $2$,
in which the successive images   of
a (very well chosen) open rectangle have a diameter that converges to $0$,
contradicting the naive intuition that 
the expanding property leads to
evergrowing iterates and even to exactness.

In view of the  shape and nonlinearity of the recurrence relation, the technics developed in \cite{SAU} seemed perfectly adapted. In this work,
the author studies expanding, piecewise defined
transformations $T$ of a compact subset $\Om$ of  $\R^k$.
\\
Under 
geometric conditions, he  proves the existence of an absolutely continuous
invariant measure $\mu$ (ACIM) for $T$. Moreover, if $T$ is mixing,
the  correlations decay at an exponential rate. That means that, for
 convenient functions $f,h$,  there exist $\Lambda \in [0,1[$ and a 
 real constant $C(f,h)$ such that, for any $n$,
 \[
\left \vert \int_{\Omega}f\circ T^n h \ d\mu -
\int_{\Omega}f\ d\mu\int_{\Omega} h \ d\mu
 \right \vert\leq C(f,h) \Lambda^n.
  \]
  In other words, if $Y_0 \sim \mu$ and if $Y_n=T^nY_0$, 
  $ {\rm cov}(f(Y_n),h(Y_0))\leq  C(f,h) \Lambda^n$.
  
Obtaining the mixing property is non trivial. The classical
conditions on the
transfer operator
or the more recent random covering condition \cite{ANV} are difficult to
obtain. For example,
in \cite{SOC} the author even obtains exactness, but the shape of the partition
is very regular. Therefore, the mixing condition has, up to now, to remain an
additional condition.

To apply the results of \cite{SAU} in our situation,
 one needs to  derive  a 
 transformation $T$ of a compact subset of  $\R^k$ from the function 
  $\Phi_0$ leading the process \eqref{pb1}.
 This was achieved, in the deterministic case,
 in \cite{JMN1,JMN2} for recurrence relations of all orders $k\geq 2$. 

The present work is the following step.
It deals with the full problem
\eqref{pb1}, with its stochastic perturbation $\Theta$.
Instead of treating $\Theta$ as a perturbation, we study
 a family $(\ph_{\theta})_{\theta\in \R}= (\Phi_0  + \theta)_{\theta\in \R}$
 of functions leading the process and 
 we consider that $\theta$ is a realization of the random variable  $\Theta$,
 with (for example) normal distribution. Each $\ph_{\theta}$ gives rise 
 to a transformation $T_{\theta}$ - to which we apply the results of \cite{SAU}- and to its transfer operator $ P_{\theta}$. \\
 Following the lines of \cite{ANV},  we
 then study an {\it averaged transfer operator} $P$ which is,
 in some sense, the expected value of 
 all the individual transfer operators  $ P_{\theta}$ (under the
 law of the perturbation $\Theta$). If the
individual transformations $T_{\theta}$ satisfy the conditions in \cite{SAU}
 uniformly in $\theta$, the operator $P$ then satisfies the conditions of \cite{ITM}. As a consequence, it
has a spectral decomposition, a spectral gap, all properties which ensure
the existence of an invariant measure in a certain sense.
This corresponds to  the annealed   point of view
 and we are in the conditions of Example 2.7 of  \cite{ANV}.
 The averaged operator is defined in the same way as the Foias
operator in \cite{LM}, section 12. But, since the transformation lacks
a continuity property (for a fixed value of the perturbation, our
transformation is piecewise $C^2$ and not continuous), the conclusions
valid for the Foias operator do not hold in our case.

The article is organized as follows.

Section \ref{sec-hyp-res} gives the precise setting and hypotheses. 
The functions  $\ph_{\theta}$ are all additive perturbations of the same function 
 $\Phi_0 $  defined on a compact hypercube
 $[-L,L]^k$ of $\R^k$. Since each $\ph_{\theta}$ is  supposed to take values in $[-L,L]$, 
it is piecewise defined on  $[-L,L]^k$ and the precise definition domains
 depend on $\theta$.
We then define the 
compact set $\Om$ and the
transformations $T_{\theta}$ derived from the $\ph_{\theta}$ and to which we apply
 the results of \cite{SAU}. The  $T_{\theta}$ and the set  are defined using a
small parameter $\gamma$, which does not depend on $\theta$. To ensure its
 existence, we need  $\Phi_0$ and its derivatives to satisfy 
the precise, quantitative conditions \eqref{cond}.
We then give the main  results. Theorem \ref{unif} states that the conditions of \cite{SAU} are uniformly satisfied, Theorem \ref{av-TPF} gives the decomposition of the averaged transfer operator $P$ and Theorem \ref{th-inv-meas} gives the invariant measure for the $k$ dimensional process and for the original $\R$ valued process.  
The proofs are detailed in the following sections.

In Section \ref{sec-unif},
we use the same technics as in \cite{JMN1, JMN2} to construct
the individual transformations $T_{\theta}$ of a compact set $\Omega$  of $\R^k$
 and thus to
imbed the recurrence relation into 
$\R^k$. Since the transformations depend on $\theta$, 
   we have to check carefully that 
the conditions stated in \cite{SAU} are satisfied uniformly with respect to
$\theta$. This is particularly delicate when it comes to the precise
definition domains of $\ph_{\theta}$. The results of this section are
summarized in Theorem \ref{unif}.

In Section \ref{sec-avPF},
according to \cite{ANV} 
we define and study an averaged transfer operator $P$. 
The individual transfer operators $P_{\theta}$ 
have uniformity properties, for example Section \ref{sec-unif} proves that 
they satisfy a Lasota-Yorke inequality 
with constants independent of $\theta$. The same inequality holds for $P$.
Applying the result of \cite{ITM}  to the averaged $P$ gives a spectral decomposition
 of $P $ as the
 sum of projections on finite dimensional spaces and of a ``small'' residual operator $S$ (small in that it has a small 
spectral radius). The spectrum of $P $ then consists of a finite number of 
modulus $1$ eigenvalues, all of them simple, and of the spectrum of 
$S$, which is contained in a ball of radius $<1$, centered on the origin. This particularity is called the spectral gap.

Denote by $\Pr_{\Theta}$ the tensor product
of the probability measures $\Pr_{\Theta_0} $  on $\R^{\otimes \N}$.
One of the projections is associated with the eigenvalue $1$.
Applying it to the (normalized) constant function on $\Omega$ gives an
invariant measure $\mu$, in the sense that the measure
$\Pr_\Theta \otimes \mu$ is invariant for the mapping
\begin{equation}\label{Effe}
\begin{matrix}
 F : & \R^\N \times \Omega & \longrightarrow & \R^\N \times \Omega \\
& (\vartheta, x) & \longmapsto & \left( \mathcal S \vartheta, T_{\theta_0} (x)
\right)
\end{matrix}
\end{equation}
where $\mathcal S$ denotes the left shift. It remains to compute the marginal
distributions to get an invariant measure for the original $1$-dimensional
process.


\section{Hypotheses and main results}\label{sec-hyp-res}
\subsection{Notations, concepts, hypotheses}\label{sec-nch}

We now present the problem more precisely.
Let $L>0 $, let 
the stochastic process $(X_n)$, bounded and with values
 in $ [-L,L]$,  be defined as follows. 
Let $\beta$ be a  positive real number and  ${\Phi_0}$ be a real-valued,
  $C^2$ application,
 defined on a neighbourhood $(-L-\beta,L+\beta)^k$
of $[-L,L]^k$.\\
 We  introduce a family  $\{ \ph_{\theta}\}_{\theta \in \R}$ of applications
 indexed by $\theta\in \R$. 
Each $\ph_{\theta}$ is defined piecewise on  $[-L,L]^k$ as an additive perturbation
 of $\Phi_0$, with values in 
 $[-L,L]$. 
More precisely, let 
$(O^{\theta}_j)_{j\in \Z,\theta\in \R}$ be a family of open subsets of $ (-L,L)^k$
which split  $[-L,L]^k$ according to the values taken by  $\Phi_0$ and
denote by $S^{\theta}_j $ the parts of their boundaries which are inside
$ [-L,L]^k$ :

\begin{equation}\label{OjSj}
\begin{array}{lll}\displaystyle
O^{\theta}_j=\{ x\in (-L,L)^k,\  (2j-1)L < \theta + {\Phi_0}(x) < (2j+1) L \}
\\ \\
\displaystyle
S^{\theta}_j  =\{ x\in [-L,L]^k,\  (2j-1)L = \theta + {\Phi_0}(x) \}.
\end{array}
\end{equation}

For every $\theta\in \R$ and every $x\in [-L,L]^k$, set
\begin{equation}\label{def-phi-om}
\begin{array}{lllll}
\ph_{j,\theta}(x) & = & {\Phi_0}(x)+ \theta \quad [-L,L]\\
   &  = & {\Phi_0}(x)+ \theta -2jL, & j=
 \lfloor \frac{{\Phi_0}(x)+\theta +L}{2L} \rfloor,
\end{array}
\end{equation}
where $ \lfloor x\rfloor$ denotes the integer part of $x$.
We define $\ph_{\theta}$ on $[-L,L]^k$ (up to a null set) by
\begin{equation}\label{phi-om}
\ph_{\theta}|_{O^{\theta}_j}= \ph_{j,\theta},
\end{equation}
provided $O_j^{\theta}$ is not empty. 

If, at time $n\geq k$, the perturbation $\Theta_n$
takes the real value $\theta$, 
the next state $X_{n+1}$ is then given by
\begin{equation}\label{pbs}
X_{n+1 }=\ph_{\theta}(X_{n-k+1},\dots, X_{n-k}).
\end{equation}
We imbed the recurrence relation \eqref{pbs}
in $\R^k$, setting 
\begin{equation}\label{imbed}
Y_n=(X_{n-k+1},\dots , \gamma^{j-1} X_{n-k+j},\dots \gamma^{k-1}X_n),
\end{equation}
where $\gamma$ is a small positive parameter, to be chosen.

The conditions \eqref{cond} listed below ensure that there exists a 
 $\gamma$, independent on  $\theta\in \R$ and small enough to allow
us to apply  Theorem 5.1 in \cite{SAU}.  The natural way, with
 $\gamma=1$, doesn't allow the transformation $T_{\theta}$ to
be expanding (condition (PE4) of Theorem \ref{unif}). Anyway,
  condition (PE5) on a geometric constant $\eta$ (see \eqref{eta}) imposes
 even stricter conditions. 

Let us first  define the applications and their domains. We give the convenient
 value of $\gamma$ in \eqref{cond}.

One sets
\begin{equation}\label{Omega}
\Omega= \prod_{j=1}^k [-\gamma^{j-1}L, \gamma^{j-1}L]. 
 \end{equation}

One denotes by $\Gamma$ the dilation linking $\Omega$ and $[-L,L]^k$ :
\begin{equation}\label{Gamma-dilate}
\begin{array}{llll}
\Gamma : & \R^k & \rightarrow &\R^k \\
  & x=(x_j)& \mapsto & (x_j \gamma^{-j+1}), \\
\end{array}
\end{equation}
hence $\Omega =  \Gamma^{-1}([-L,L]^k)$.

One then sets 
\begin{equation}\label{ouv-front}
\begin{array}{llll}
\displaystyle \forall r\in [0,\beta[,\quad
\Omega_{r} = \Gamma^{-1}((-L-r,L+r)^k) =
  \prod_{s=1}^k (-\gamma^{s-1}(L+r),  \gamma^{s-1}(L+r)L) ,\\
\displaystyle
U^{\theta}_j= \Gamma^{-1}(O^{\theta}_j)   =\left\{ 
x\in \Omega_0,\ (2j-1)L<\theta + \Phi_0(\Gamma(x))< (2j+1)L   \right\}\\ \\
\displaystyle
\Sigma^{\theta}_j  =\Gamma^{-1}(S^{\theta}_j )  =\left\{  x \in \Omega,\ 
 (2j-1)L = \theta + \Phi_0(\Gamma(x)) \right  \}  
\end{array}
\end{equation}

The open sets $\Omega_{r}$ are neighbourhoods of $\Omega$ and appear in
Section \ref{sec-unif}.
The open set $U_j^{\theta}$  is empty if and only if 
$O_j^{\theta}$ is empty. 
The properties of $\Omega$, of the $U^{\theta}_j$ and of the slanted parts of 
the boundaries, the $\Sigma^{\theta}_j$ 
 (regularity, number) do not depend on the value of $\gamma >0.$
\\ 

On the set $\Omega_{\beta}\setminus \{ x \in\Omega_{\beta}, 
\theta + \Phi_0(\Gamma(x)) \in (2\N+1)L
 \} $,
the transformation $T_{\theta}$ is defined  by
\begin{equation}\label{Tomega}
T_{\theta}(u) = \left( \frac{u_2}{\gamma}, ~ \dots, \frac{u_k}{\gamma},
 ~ \gamma^{k-1} \ph_{\theta}(u_1, \frac{u_2}{\gamma}, ~ \dots, ~ \frac{u_k}{\gamma^{k-1}})\right),
\end{equation}
with $\ph_{\theta}$ given by \eqref{phi-om}. 

We introduce an auxiliary function $S$ on  $\Omega_{\beta}$ :
\begin{equation}\label{Esse}
S(u) = \left( \frac{u_2}{\gamma}, ~ \dots, \frac{u_k}{\gamma},
 ~ \gamma^{k-1}
 \Phi_0(u_1, \frac{u_2}{\gamma}, ~ \dots, ~ \frac{u_k}{\gamma^{k-1}})\right)=
\left( \frac{u_2}{\gamma}, ~ \dots, \frac{u_k}{\gamma},
 ~ \gamma^{k-1}
 \Phi_0(\Gamma(u))\right).
\end{equation}

We then define   $T_{j}^{\theta}$ on  $\Omega_{\beta}$ by
\begin{equation}\label{T-et-S}
T_{j}^{\theta}(u) = S(u)+ (0,\dots, 0, \gamma^{k-1}(\theta -2jL)).
\end{equation}
 One sees that  
$\left. T_\theta \right|_{U_j^\theta} = \left. T_{j}^{\theta} \right|_{U_j^\theta}$.

From now on, we are interested in the $k$-dimensional valued process
defined by its first term $Y_0$ and the relation
$$
Y_{n+1}= T_{\theta}(Y_n),
$$
where the relationship between the $\R^k$ valued processus $Y=(Y_n)$ and
the real valued $X=(X_n)$
is given by \eqref{imbed}:
\begin{equation*}
Y_n=(X_{n-k+1},\dots , \gamma^{j-1} X_{n-k+j},\dots \gamma^{k-1}X_n).
\end{equation*}

Lemma \ref{S-diff} below proves that  
 $T_{j}^{\theta}$ is a global diffeomorphism on
 $\Omega_{\beta}$ onto its image  $T_{j}^{\theta}(\Omega_{\beta})$. Hence, 
 $T_{j}^{\theta}$  is a diffeomorphism  on 
 $U^{\theta}_j$ onto  $T_{j}^{\theta}(U^{\theta}_j)$. This holds for 
any bounded open subset of 
  $\Omega_{\beta}$.
Later on, we will need to restrict  $T_{j}^{\theta}$ to more specific open sets
(defined by \eqref{Vjomega}),
in order  to apply directly 
 the results of \cite{SAU}.\\

\begin{figure}[h]
\centering
\includegraphics[scale=0.3]{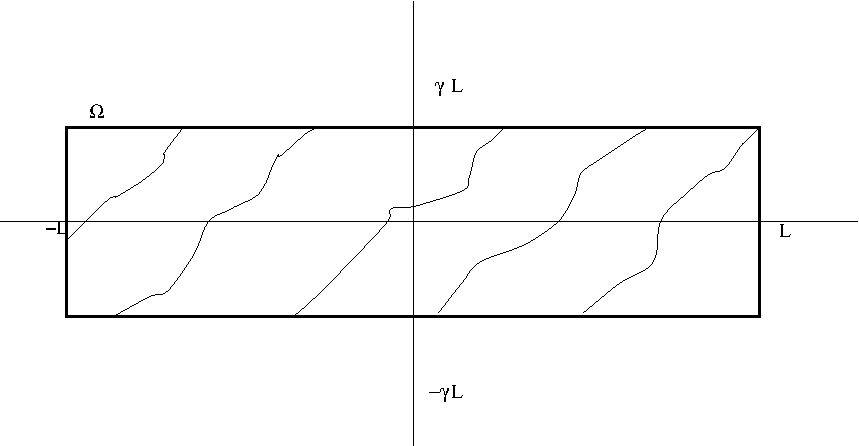}
\caption{\label{cross} Open sets $U_j^{\theta}$ with boundaries and crossings.
}
\end{figure}

A central parameter in \cite{SAU} is the {\it number of crossings}, that is,
the number  of regular parts of the
 boundaries of the open sets of $\Omega$ 
which intersect at one point. 
In the present case, the boundaries consist of parts of hyperplanes
(the boundaries of
 $\Omega$ itself) and of the slanted parts $\Sigma_j^{\theta}$. 
Whatever the value of the parameter $\gamma$, the  $\Sigma_j^{\theta}$ are
 regular under the conditions \eqref{cond} (see 
 Lemma \ref{front-neg}).  
 There are at most $k$ hyperplanes
  crossing at one point, and one slanted part (see Figure \ref{cross} in the 
case when $k=2$). Since the 
 $\Sigma_j^{\theta}$ are common to 
two open sets, they count double (see Lemma 2.1 \cite{SAU}), which
 proves that the maximal number of crossings is
$k+2$ for a countable number of values $\theta$  and $k+1$ 
for the rest of the values.  \\

Let $Y$ be the maximal number of crossings. For a real number 
 $\sigma>1$, we define
\begin{equation}\label{eta}
\eta_{0,Y}(\sigma)=
\eta_0(\sigma, Y) = \frac{1}{\sqrt{\sigma}} +
\frac{4}{\sqrt{\sigma} - 1} Y \frac{\gamma_{k-1}}{\gamma_k},
\end{equation}
where $\gamma_n$ is the volume of the $n$-dimensional unit ball in $\R^n$.
Let  $S_k(Y)$ be the real number for which 
$\eta_0(\sigma, Y)=1$:
\begin{equation}\label{seuil}
S_k(Y) = \eta_{0,Y}^{-1}(1).
\end{equation}

We choose the dilation coefficient  $\sigma$:
\begin{equation}\label{sigma}
\sigma > S_k(Y)>1
\end{equation}
and  two real constants  $C_1,C_2>1$. We now  prescribe the conditions: 
\begin{equation}\label{cond}
\left\{
\begin{array}{llllll}
&
\displaystyle \gamma^{-2} & =& C_1 \sigma, \\ \\
\displaystyle
\forall x \in (-L-\beta, L+\beta)^k ,&
\displaystyle
\left(\frac{\partial {\Phi_0}}{\partial x_1}(x)\right)^2&
 \geq& C_1^{k-1}C_2\sigma^k ,\\
 \\ 
\displaystyle
\forall x \in (-L-\beta, L+\beta)^k ,\forall j\geq 2, &
\displaystyle
\left(\frac{\partial {\Phi_0}}{\partial x_j}(x) \right)^2 &\leq 
& \displaystyle  \frac{(C_1-1)(C_2-1)}{k-1}\sigma.
\end{array}\right.
\end{equation}

This gives the value of the parameter $\gamma$, which, in turn,
defines the domain 
$\Omega$ and the transformations. 

The above conditions prove that each transformation \eqref{Tomega} 
is non singular in $L^1(\Omega, m)$, which implies the existence of
a Perron-Frobenius operator $P_{\theta}$. To study the transfer operators
$P_{\theta}$ or the averaged operator defined in \eqref{PF_moyen}, we use the
central result of \cite{ITM}. Thus, we need to introduce
 the following space $V_1$. We just recall here its definition. 
Its properties are studied in \cite{SAU}.

For every non-null Borel set $S$ of
 $\R^k$, for every $f \in L^1_m(\R^k,\R)$, one sets :
$$ \osc(f,S) = \supess_S f - \infess_S f $$
where $\supess_S$ and $\infess_S$ are the essential supremum
 and infimum on $S$ with respect to the Lebesgue measure $m$.\\
One then defines the semi norm  $| \cdot |_{1}$ by
\begin{equation}\label{snorm1}
|f|_{1}= \sup_{0<\eps<\eps_1}\eps^{-1} \int_{\R^k} 
 {\rm Osc}(f,B(x,\eps))\ dx 
\end{equation}

and the set $V_1$ by
\begin{equation}\label{V1}
V_{1} (\Omega) = \{ f \in L^1_m(\R^k,\R), ~  \supp f \subset \Omega \text{ and } \abs{f}_{1} < +\infty \}.
\end{equation}

Then, for $f\in V_1$,  the norm $ \| \  \|_{V_1}$ is defined by
$$ 
 \| f \|_{V_1} = \| f \|_{L^1_m} + |f|_{1}.
$$


\subsection{Main results}\label{ssec-result}
Let $m$ denote the Lebesgue measure on $\R^k$ and its subsets. For 
$X \subset \R^k$ and $\eps > 0$, we denote by  $B_{\eps}(X)$ the following
open set
$B_{\eps}(X) \coloneqq \left\{ x \in \R^k : d(x, X) < \eps \right\}$.

\begin{theo}\label{unif}
Under the conditions of the preceding paragraph, there exist positive
 real numbers $r_0,\eps_m$ and $K$ and a sufficiently small 
 $\eps_0 >0$ such that, for every 
 $\theta\in \R$, the following conditions hold:
\begin{itemize}
\item[{\bf (PRE)}]
There exist a finite family $(U_j^{\theta})_j$ of pairwise disjoint open subsets
 of $\Om$, a finite family $(V_j^{\theta})_j= B_{\eps_m}(U_j^{\theta})$  
such that 
$\overline{ U_j^{\theta}}\subset  V_j^{\theta}$ and $ C^2$ applications   
 $T_j^{\theta} :  V_j^{\theta} \rightarrow \R^k$, satisfying all the following 
conditions:
\item[{\bf (PE1)}]
$\forall j$,
 $ \left.T^{\theta}\right|_{U_j^{\theta}} =  \left.T_j^{\theta}\right|_{U_j^{\theta}} $ and
$B_{r_0}(T^{\theta}(U_j^{\theta})) \subset T_j^{\theta}(V_j^{\theta})$ ;
\item[{\bf (PE2)}]
$\forall j, T_j^{\theta}$ is a $C^2$ diffeomorphism on $V_j^{\theta}$ onto its image. 
Moreover, 
the determinant of the Jacobian matrix
 is uniformly Lipschitz continuous, in the following sense:\\
$
\forall j,\forall \eps\in (0,\eps_0),$ $ \forall z\in T_j^{\theta}(V_j^{\theta}), \forall x,y\in B(z,\eps)\cap T_j^{\theta}( V_j^{\theta}),
$
$$
\left|  \det\left(D(T_j^{\theta})^{-1}(x)\right) -
 \det\left(D(T_j^{\theta})^{-1}(y)  \right)\right|
\leq K \left|  \det\left(D(T_j^{\theta})^{-1}(z)\right) \right|\eps\  ;
$$
\item[{\bf (PE3)}]
$ m(\Omega\setminus \bigcup_j U_j^{\theta}) =0$;
\item[{\bf (PE4)}] 
for all  $u,v$ in $T_j^{\theta}(V_j^{\theta})$ such that $d(u,v)\leq \eps_0$,
$$
d\left( (T_j^{\theta})^{-1}(u),  (T_j^{\theta})^{-1}(v) \right)
\leq \frac{1}{\sqrt{\sigma}} d(u,v),
$$
with $\sigma >1$ defined by \eqref{sigma}. 
 The map $T_j^{\theta}$ is thus expanding;
\item[{\bf (PE5)}]
Set $G(\eps,\eps_0)= \sup_x G(x,\eps,\eps_0)$  with
$$
G(x,\eps,\eps_0)= \sum_j
 \frac{ m\left( 
 (T_j^{\theta})^{-1}\left[ B_{\eps}(\partial T_j^{\theta}(U_j^{\theta})) \right]
\bigcap 
B\left(x,\left[ 1-\frac{1}{\sqrt{\sigma}} \right]\eps_0\right)
\right) }
{m\left(B\left(x,\left[ 1-\frac{1}{\sqrt{\sigma}} \right]\eps_0\right)\right)}.
$$
Then $\eta$ defined by
\begin{equation}\label{eta-G}
\eta(\eps_0)= \frac{1}{\sqrt{\sigma}} +
 2 \sup_{\eps\leq \eps_0}\frac{G(\eps,\eps_0)}{\eps}\eps_0
\end{equation}
satisfies $\sup_{\delta\leq \eps_0}\eta(\delta) <1$.
\end{itemize}
\end{theo}

To sum up, the conditions of Theorem $5.1$ in \cite{SAU} are satisfied,
with constants $r_0,\eps_m, K,\eps_0$ independent of  $\theta\in \R$. 
In particular, condition  (PE2) proves that the transformations $T_{\theta}$ 
are non singular. For each real number $\theta$, we may define  the
 Perron-Frobenius operator associated with $T_{\theta}$ by 
$$
\forall A\in \bbb(\Omega),\ \forall f \in L^1_m(\Omega),
 \quad  \int_{T_{\theta}^{-1}(A)} f \, \d m = \int_A P_{\theta} f \, \d m.
$$
One may check that $P_{\theta} f$ has the following expression:

\begin{equation}\label{individual-PF}
 P_{\theta} f(y)= \sum_{j\in \Z, U_j^{\theta}\neq \emptyset} 
\1_{T_{j}^{\theta}(U_j^{\theta})}(y) f((T_{j}^{\theta})^{-1}(y))
 \left| {\rm det}( D(T_{j}^{\theta})^{-1}(y) ) \right|.
\end{equation}

According to Lemma 4.1 in \cite{SAU}, if $\eps_0$ is small enough, there exist 
constants $\eta\in [0,1[$ and $D\in \R^+$ such that, for all $f\in V_1$ and all
$\theta \in \R$, 
\begin{equation}\label{ineq-LY}
|P_{\theta}f|_1 \leq \eta |f|_1  + D \|f\|_{L^1_m}.
\end{equation}
Since these constants are expressed explicitely
 in terms of $\eta(\eps_0)$, $\sigma$, $K$, they do not depend on $\theta$.

We formally define  the averaged Perron-Frobenius operator associated with
this family of transformations by
\begin{equation*}\label{av-PF}
P = \int_{\R} P_{\theta} \, \d \Pr_{\Theta}(\theta) = \int_{\R} P_{\theta} \ H(\theta) \, \d\theta,
\end{equation*}
where $\Pr_{\Theta}$ is the distribution of the perturbation $\Theta$ and 
$H$, its density function. The properties of this transfer operator are developed in 
  \ref{sec-avPF}, where a proper definition is given.

This allows one to state the following result 
\begin{theo}\label{av-TPF}
Under the conditions \eqref{eta} to \eqref{cond} and with the notations
 of the preceding paragraph, 
\begin{enumerate}
\item\label{eigs}
$P$ has finitely many eigenvalues $\lambda_1=1, ..., \lambda_p \in \C$ of
 modulus $1$. 
 The corresponding eigenfunctions belong to $V_1$.
\item \label{projs}
There exist $p$ finite-rank operators $\Pi_1, ..., \Pi_p \in \lin(V_1)$, as well as $S \in \lin(V_1)$, such that
\[
P = S + \sum_{j = 1}^p \lambda_j \Pi_j .
\]
Furthermore, for all distinct $i, j \in \{1, .., p\}$, $\Pi_i^2 = \Pi_i$, $\Pi_i \Pi_j = 0$ and $\Pi_i S = S \Pi_i = 0$.
\item\label{spec-gap}
 There exist $M > 0$ and $h > 0$ satisfying 
\[
\forall n \ge 1, \; \norm{S^n}_{\lin(V_1)} \le \frac M{(1+h)^n}.
\]
\end{enumerate}
\end{theo}

\begin{theo}\label{th-inv-meas}
Recall what we denote by $F$:
\begin{equation*}
\begin{matrix}
 F : & \R^\N \times \Omega & \longrightarrow & \R^\N \times \Omega \\
& (\vartheta, x) & \longmapsto & \left( \mathcal S \vartheta, T_{\theta_0} (x)
\right)
\end{matrix}
\end{equation*} where $\Theta=(\Theta_n)_{n\in \N}$ is a sequence of i.i.d.
perturbations and $\mathcal S$ is the left shift.
Denote by $\Pr_{\Theta}$ the tensor product $ \Pr_{\Theta_0}^{\otimes \N}$.
Let $h_* = \frac{1}{m(\Omega)}\Pi_1 \1$, with $\Pi_1$ the rank one
operator associated with $\lambda_1=1$. \\
In the sense of \cite{ANV}, 
the averaged transfer operator $P$ admits an absolutely continuous invariant measure or stationary measure $\mu$, 
which has density $h_*$ with respect to the (normalized) Lebesgue measure on
$\Omega$. This means that the measure $\Pr_{\Theta}\otimes \mu$ is
invariant by the function $F$ recalled above. \\
Then, if $(\check{\vartheta}, Y)$ has distribution $\Pr_\Theta\otimes \mu $ and if
the  $\R^k$-valued processus $Y=(Y_n)$ and
the real valued $X=(X_n)$ are related as described in \eqref{imbed},
$X_n$ has for density any one  of the following densities
\begin{equation}\label{dens-for-X}
t\mapsto \gamma^{j-1}
\int_{\prod_{i\leq k,i\neq j}[-\gamma^{i-1}L,\gamma^{i-1}L] }
h_* \left(x_1, ... , x_{j-1}, \gamma^{i-1}t, x_{j+1}, ..., x_k \right) \, \d x_1 ... \d x_{j-1} \d x_{j+1} ... \d x_k
\end{equation}

\end{theo}

Let us prove the last fact. 
If 
$(\check{\vartheta}, Y)$ has distribution $\Pr_\Theta\otimes \mu $,
 so has $F(\check{\vartheta}, Y)$. If the initial term
$(\check{\vartheta}_0, Y_0)$  has distribution  $\Pr_\Theta\otimes \mu $,
so have all terms $F^n(\check{\vartheta}_0, Y_0)= (S^n\check{\vartheta}_0, Y_n)$

This implies that $Y_n$ has the density $h_*$. Considering  the marginal
distributions we obtain that
$X_{n-k+1}$ admits, as an density,  the first marginal density of $h_*$,
$\gamma ^{j-1} X_{n-k+j}$ admits the $j$-th  marginal density of $h_*$ for $j=1, \dots, k$. This gives \eqref{dens-for-X}, exactly as in 
Theorem 2 of \cite{JMN2}.

\section{Proof of Theorem \ref{unif}}\label{sec-unif}

\subsection{Expanding property}\label{section_dilatance}


In this part we prove that, under the conditions listed in the preceding
section, the application $S$ defined by \eqref{Esse} is expanding
on $\Omega_{\beta}$. Since they differ from $S$ by an additive constant, the 
applications $T_{j}^{\theta}$ share this property.
We first express the difference $S(u)-S(u')$ thanks to differential calculus. 
There appears a matrix $\B$, the eigenvalues of which measure the expansiveness.
We prove that the eigenvalues of $\B$ are greater than or equal to $\sigma$,
 which gives the inequality of Proposition \ref{dilatance}.

Note a difference with  \cite{JMN2}.  Here we compute explicitly the 
spectrum of $\B$, whereas we only located it in the
Gerschg\"orin disks. 
As a consequence,  the conditions \eqref{cond}
required  here are simpler than in \cite{JMN2}: they consist in estimates
 about the first derivatives of
$\Phi_0$, instead of estimates about products of the 
same derivatives.

\begin{prop}\label{dilatance}
Suppose that  $\Phi_0$ and  $\gamma$ satisfy  the conditions \eqref{eta}
to \eqref{cond}.
 Then the function $S$ defined by \eqref{Esse} satisfies
\begin{equation*}\label{wished-ineq} 
\forall (u,u')\in  \Omega_\beta^2,\quad
\norm{S(u) - S(u')}^2 \ge \sigma \norm{u - u'}^2.
\end{equation*}
\end{prop}

{\bf Proof.}\\
Let $u, u'$ belong to  $ \Omega_\beta$. Then 
\[
\norm{S(u) - S(u')}^2 = \sum_{i = 2}^k \gamma^{-2} (u_i - u_i')^2 + \gamma^{2k-2}
 \left( (\Phi_0\circ\Gamma)(u) - (\Phi_0\circ\Gamma)(u') \right)^2.
\]
There exists $t_0 \in [0, 1]$ such that
\begin{align*}
(\Phi_0\circ\Gamma)(u) - (\Phi_0\circ\Gamma)(u') &= \left.\dfrac{\d}{\d t}
 (\Phi_0 \circ \Gamma)\left(tu + (1-t)u'\right)\right|_{t=t_0} \\
&= \Gamma(u - u') \cdot \nabla (\Phi_0 \circ \Gamma) (w) = \sum_{i = 1}^k 
\gamma^{-i+1} (u_i - u_i') \dfrac{\partial \Phi_0}{\partial x_i}(\Gamma w)\ ,
\end{align*}
where we set $w = t_0 u + (1 - t_0) u'$. We obtain 
\begin{align*}
\norm{S(u) - S(u')}^2 &= \sum_{i = 2}^k \gamma^{-2} (u_i - u_i')^2 + \gamma^{2k-2}
 \sum_{i = 1}^k \sum_{j = 1}^k \gamma^{-i-j+2} (u_i - u_i')(u_j - u_j') 
\dfrac{\partial \Phi_0}{\partial x_i}(\Gamma w)
 \dfrac{\partial \Phi_0}{\partial x_j}(\Gamma w) \\
&= \prescript t{}{(u - u')} \B (u - u')
\end{align*}
where the matrix $\B = \B(\Gamma w)$ is given by

\[
\B (x) = \frac1{\gamma^2} \begin{pmatrix}
0 &  \\
& I_{k-1}
\end{pmatrix} + \begin{pmatrix}
\gamma^{k-1} \dfrac{\partial \Phi_0}{\partial x_1} (x) \\
\gamma^{k-2} \dfrac{\partial \Phi_0}{\partial x_2} (x) \\
\vdots \\
\dfrac{\partial \Phi_0}{\partial x_k} (x)
\end{pmatrix} \begin{pmatrix}
\gamma^{k-1} \dfrac{\partial \Phi_0}{\partial x_1} (x) & \gamma^{k-2} \dfrac{\partial \Phi_0}{\partial x_2} (x) & \cdots & \dfrac{\partial \Phi_0}{\partial x_k} (x)
\end{pmatrix}.
\]
Set
\begin{equation}\label{v0V}
v_1(x) \coloneqq \gamma^{k-1} \dfrac{\partial \Phi_0}{\partial x_1} (x), \hspace{1em}
V(x) \coloneqq \prescript t{}{\begin{pmatrix}
	\gamma^{k-2} \dfrac{\partial \Phi_0}{\partial x_2} (x) & \cdots & 
	\gamma \dfrac{\partial \Phi_0}{\partial x_{k-1}} (x) &
	\dfrac{\partial \Phi_0}{\partial x_k} (x)
\end{pmatrix}}\in \R^{k-1}.
\end{equation}

The proof of the proposition is then a consequence of the following result.
\begin{prop}\label{eig-B}
For all $x\in (-L-\beta, L+\beta)^k$, the eigenvalues $\lambda$ of $B=B(x)$ 
satisfy
$
\lambda \geq \sigma.
$
\end{prop}

{\bf Proof of Proposition \ref{eig-B}.}\\
Unsing the notations \ref{v0V}, $B$ can be written as 
\[
{\B} = \gamma^{-2} \begin{pmatrix}
0 & 0 \\
0 & I_{k-1}
\end{pmatrix} + \begin{pmatrix}
v_1^2 & v_1 \prescript t{}V \\
v_1 V & V \prescript t{}V
\end{pmatrix}.
\]
 The functions $v_1$ and
 $V$ are taken in $x=\Gamma w$. \\

When $V = 0$, $\B$ is  diagonal with two eigenvalues: a simple one,
$v_1^2$ and a multiple one,  $ \gamma^{-2}$. The conditions \eqref{cond}
imply that they are $>\sigma$, since 
$ v_1^2(x) \ge C_2 \sigma $ and $\gamma^{-2} = C_1\sigma$.

We now treat the case when $V\neq 0$.
  
\begin{lem}\label{matrice}
Let $x \in \left(-L-\beta, L+\beta\right)^k$ be such that
$V \coloneqq V(x) \ne 0$. Then:
\begin{enumerate}
\item $\gamma^{-2}$ is an eigenvalue of $\B$, with multiplicity 
 $k - 2$. Its eigenspace is the set of vectors 
  $\begin{pmatrix}0 \\ X\end{pmatrix} \in \R^k$ with $X  $ in the orthogonal of  $ \{V\}$ (in $\R^{k-1}$).
\item The real numbers
\[
\lambda_- \coloneqq \frac{\gamma^{-2} + v_1^2 + \norm{V}^2 - \sqrt{\left(\gamma^{-2} + v_1^2 + \norm{V}^2\right)^2 - 4v_1^2\gamma^{-2}}}2,
\]
\[
\lambda_+ \coloneqq \frac{\gamma^{-2} + v_1^2 + \norm{V}^2 + \sqrt{\left(\gamma^{-2} + v_1^2 + \norm{V}^2\right)^2 - 4v_1^2\gamma^{-2}}}2,
\]
are simple eigenvalues, associated with the eigenvectors
 $\begin{pmatrix}t_\pm\\V\end{pmatrix}$, where $t_\pm = (1-\lambda_\pm^{-1}\gamma^{-2}) v_1$.
\end{enumerate}
\end{lem}

{\bf Proof of the lemma.}\\
\noindent{\it Step 1.}\
For every vector  $\begin{pmatrix}t \\ X\end{pmatrix} \in \R^k$, one
checks that
$$
\B \begin{pmatrix}t \\ X\end{pmatrix} = \gamma^{-2} \begin{pmatrix}t \\
 X\end{pmatrix} \qquad \iff
\quad(v_1 t + \prescript t{}V X)\begin{pmatrix}v_1 \\
 V\end{pmatrix}  = \gamma^{-2} \begin{pmatrix}t \\
 0\end{pmatrix} \ . \label{vp_gamma-2}
$$
This implies that
 $v_1 t + \prescript t{}V X = 0$ since $V \ne 0$. Hence $t = 0$ by \eqref{vp_gamma-2}, which gives  $v_1 t + \prescript t{}V X = \prescript t{}V X = 0$ and $X \in \{V\}^\bot$. Reciprocally, every vector
 $\begin{pmatrix}0 \\ X\end{pmatrix}$ with $X \bot V$ satisfies \eqref{vp_gamma-2}. This proves that $\gamma^{-2}$ is an eigenvalue of $\B$, 
 with multiplicity   $k - 2$, associated with the eigenspace
\[
E(\gamma^{-2}) = \left\{ \begin{pmatrix}0 \\ X\end{pmatrix} : X \in \R^{k - 1}, \ X \bot V \right\}.
\]
\noindent{\it Step 2.}\
We now investigate the other eigenvalues. We look for eigenvectors taking the shape $\begin{pmatrix}t\\V\end{pmatrix}$ avec $t \in \R^*$. For $\lambda \in \R, t \in \R^*$, one may write
$$
\B \begin{pmatrix}t\\V\end{pmatrix} = \lambda \begin{pmatrix}t\\V\end{pmatrix} \qquad \iff \qquad \left\{
\begin{array}{rcl}
\lambda t & = & v_1^2 t + \norm{V}^2 v_1 \\
\lambda & = & \gamma^{-2} + v_1 t + \norm{V}^2
\end{array}
\right. .$$
Remark that this system ensures that  $\lambda \ne 0$. Otherwise 
$v_1^2 t + \norm{V}^2 v_1$,  $v_1 t + \norm{V}^2 $ and finally
  $\gamma^{-2} $ would be equal to $0$, which is impossible.

We now solve the system. One sees that
$$
\left\{
\begin{array}{rcl}
\lambda t & = & v_1^2 t + \norm{V}^2 v_1 \\
\lambda & = & \gamma^{-2} + v_1 t + \norm{V}^2
\end{array}
\right.\quad \iff \quad
 \left\{
\begin{array}{rcl}
t & = & (1-\lambda^{-1}\gamma^{-2}) v_1 \\
P(\lambda) & = & 0
\end{array}
\right.
$$
with $P(X) = X^2 - (\gamma^{-2} + v_1^2 + \norm{V}^2) X + \gamma^{-2} v_1^2$.

The discriminant of $P$ is $\Delta = \left(\gamma^{-2} + v_1^2 + \norm{V}^2\right)^2 - 4v_1^2\gamma^{-2} > 0$.
 This implies that the remaining eigenvalues  are the roots of $P$,
  that is :
\[
\lambda_- \coloneqq \frac{\gamma^{-2} + v_1^2 + \norm{V}^2 - \sqrt{\Delta}}2 <
\ \lambda_+ \coloneqq \frac{\gamma^{-2} + v_1^2 + \norm{V}^2 + \sqrt{\Delta}}2
.
\]
And the two remaining eigenvectors take the shape
 $\begin{pmatrix}t_\pm\\V\end{pmatrix}$ with
 $t_\pm = (1-\lambda_\pm^{-1}\gamma^{-2}) v_1$.
This achieves the proof of Lemma \ref{matrice}
\qed\\

It remains to prove that 
$\min \{ |\lambda| : \lambda \in \operatorname{Sp}(\B) \} \ge \sigma$ in
 the case when $V\neq 0$. In this case, 
the spectrum of the matrix $\B$  consists of the eigenvalues 
$ \gamma^{-2}$ and  $\lambda_+ > \lambda_- >0$. We already know that 
 $\gamma^{-2}$ is greater than $\sigma$ and we now prove that 
$ \lambda_- \geq \sigma$. 

Set  $c(x) \coloneqq \sigma^{-1} v_1^2 (x) \ge C_2$ and  $\Delta =
 \left(\gamma^{-2} + v_1^2 + \norm{V}^2\right)^2 - 4v_1^2\gamma^{-2} > 0$. Then
\[
\begin{array}{rcl}
\lambda_- \ge \sigma &\iff& (C_1 + c)\sigma + \norm{V}^2 - \sqrt{\Delta} \ge 2
 \sigma \\
&\iff& (C_1 + c - 2)\sigma + \norm{V}^2 \ge \sqrt{\Delta} \\
&\iff& (C_1 + c - 2)^2\sigma^2 + \norm{V}^4 \\
& &+ \, 2(C_1 + c - 2)\sigma\norm{V}^2 \ge \left((C_1 + c)\sigma + 
\norm{V}^2\right)^2 - 4cC_1\sigma^2 .
\end{array}
\]
Developing  $\left((C_1 + c)\sigma + \norm{V}^2\right)^2 = (C_1 + c)^2 \sigma^2
 + 2(C_1 + c)\sigma\norm{V}^2 + \norm{V}^4$  gives the following inequality: 
\begin{align*}
&\left[(C_1+c-2)^2 - (C_1+c)^2 + 4 c C_1\right]\sigma^2 - 4\sigma\norm{V}^2
 \ge 0 \\
\iff &4 \left[(C_1 c - C_1 - c + 1)\sigma - \norm{V}^2\right] \sigma \ge 0 \\
\iff &\norm{V}^2 \le (C_1 c - C_1 - c + 1)\sigma = (C_1 - 1)(c - 1)\sigma \ .
 \tag{$**$}
\end{align*}

One may check that the conditions \eqref{cond}
and the definition \eqref{v0V} of $v_1$ and $V$ imply that, for all
 $x \in (-L-\beta, L+\beta)^k$,
\[
\hspace{1em} v_1^2(x) \ge C_2 \sigma \hspace{1em} \mathrm{ and } \hspace{1em}
 \norm{V(x)}^2 \le (C_1 - 1)(C_2 - 1) \sigma . 
\]
Hence, the condition  $(**)$ holds for an arbitrary $x=\Gamma w$ 
and we obtain that  
$\lambda_- \ge \sigma$.

We have proved that 
 $\min \{ |\lambda| : \lambda \in \operatorname{Sp}(\B) \} \ge \sigma$,
 which concludes the proof of  Proposition \ref{eig-B} and hence of 
 Proposition \ref{dilatance}.
\qed


\subsection{Number of open sets, uniformity}
This paragraph goes on with the rather long series of 
results proving  that the  conditions listed in Proposition
\ref{unif} are satisfied.\\
 It is not difficult to see that
 the number of non empty open subsets  $U_j^{\theta} $ is
bounded, see \eqref{nb-ouv}.\\
 We then state a uniform version of the implicit functions theorem. 
It is necessary because 
the constants $\eps_m,r_0$ and $\eps_0$ of Theorem \ref{unif} express
a geometric uniformity. The open sets $U_{j}^{\theta}$ 
are determined by the function $\Phi_0$, but we need  to ensure that the
slightly larger open sets  $V_{j}^{\theta}$
defined by \eqref{Vjomega} do not depend on $\theta$ in a
way we can't control. Similarly, we use open charts (see \eqref{carte})
 in a proof below
  and we need the charts at each point to contain a ball
 of the same size. 
Another consequence of this theorem is Lemma \ref{front-neg} which gives
the negligibility of the boundaries. This is 
condition (PE3) of  Theorem \ref{unif}.
This implies that the transformation $T_{\theta}$ is defined almost everywhere,
since it is defined on $\Omega \setminus \bigcup \Sigma_j^{\theta} $. 
\\

Let  $m$ and $M$ be the  minimum and the maximum of $\Phi_0$ on
$[-L,L]^k$.
Let us fix $\theta\in \R$. If the open set  $O_j^{\theta} $ is not empty, then
the index $j$ satisfies 
$$
\frac{\theta +m-L}{2L} <j<\frac{\theta +M +L}{2L}.
$$
The maximal number $N$ of non empty open sets $ O_j^{\theta} $ or 
$ U_j^{\theta} $ satisfies
\begin{equation}\label{nb-ouv}
N<\frac{M-m+2L}{2L}+1.
\end{equation}
Indeed, if  $O_j^{\theta} \neq \emptyset$, 
 then there exists $x\in [-L,L]^k$ such that
$
(2j-1)L < \Phi_0(x)+\theta < (2j+1)L.
$
This gives the inequality. 

\begin{theo}\label{TFIU}
Let $E$, $F$ and $G$ be Banach spaces. The space $E\times F$ is endowed with
the $\max $ norm. Let $O$ be an open set of 
 $E \times F $.\\
Let $f : O \rightarrow G $ be  $C^2$ at least on $O$.
Suppose that, for every 
$(a,b) \in O$, $\partial_2 f(a,b) \in Isom(F,G)$.\\
Suppose that there exist nested open sets 
 $\Omega', \Omega''$ satisfying
$$
\Omega'\subset \overline{\Omega'} \subset \Omega'' \subset \overline {\Omega ''}
\subset O,
$$
with $ Df, D^2 f$ bounded on $\overline{\Omega''}$.
Suppose that
$\left( \partial_2 f\right)^{-1}$  is bounded too on  $\overline{\Omega''}$.
Assume that the distance between 
$\overline{ \Omega'}$ and $\complement \Omega''$ is a positive real number 
$\delta >0$.\\
Then there exist 
$r>0,\rho>0$ such that  $0<\rho<r<\delta$ and satisfying, for every 
  $(a,b) \in \Omega'$, the following conditions:
\begin{enumerate}
\item There exist an open neighbourhood  $U=B_E(a,r)\times B_F(b,r) \subset
 \Omega ''$ 
of $(a,b)$ in $O$, an open neighbourhood $V$ of $a$ in $E$ and a function 
$\varphi \in C^2(V,F)$ such that, for every $(x,y) \in \Omega''$,
$$ (x,y) \in U ~ {\rm and } ~ f(x,y) = f(a,b) \Longleftrightarrow x \in V ~
 {\rm and} ~ y = \varphi(x).$$
Moreover,  $B_E(a,\rho)\subset V\subset B_E(a,r)$.
\item For all  $x \in V$, 
$$ \partial_2 f(x,\varphi(x)) \in Isom(F,G) ~ {\rm and} ~ D\varphi(x) =
 - [\partial_2 f(x,\varphi(x))]^{-1} \circ \partial_1 f(x,\varphi(x)).$$
\item If $f \in C^k (O)$, then $\varphi \in C^k (V)$.
\end{enumerate}
\end{theo}

We apply this result to the following   $C^2$ function $f$
\begin{equation}\label{f-TFI}
\begin{array}{llll}
	f : & \Omega_{\beta} &\rightarrow & \R\\ 
	& u & \mapsto & \Phi_0(\Gamma(u)).
\end{array}
\end{equation}
The global open set $O$ of Theorem \ref{TFIU} is 
$O=\Omega_{\beta}$ and the intermediate sets are 
$\Omega'= \Omega_{\beta/3}$ and $\Omega''=\Omega_{2\beta/3}$. They are nested:
$$
\Omega_{\beta/3}\subset \overline{ \Omega_{\beta/3}} \subset  \Omega_{2\beta/3}
\subset \overline{ \Omega_{2\beta/3}} \subset  \Omega_{\beta}.
$$
Denote by  $d >0$ the real number:
\begin{equation}\label{def_d}
d=\min( {\rm dist}(\complement { \Omega_{2\beta/3}},\overline{\Omega_{\beta/3}}),
{\rm dist}(\Omega,\complement \Omega_{\beta/3})),
\end{equation}
where the distance is the $\max$ distance, considering 
$\R^k= \R\times \R^{k-1}$.

We now turn to the function $f$.
Its first and second order derivatives are bounded on
 the compact set $\overline{ \Omega_{2\beta/3}}$. Moreover, since 
$\displaystyle
\frac{\partial f}{\partial u_1}(a) = 
\frac{\partial \Phi_0}{\partial u_1}(\Gamma(a)),
$
the conditions \eqref{cond} imply that  
$\displaystyle
\left| \frac{\partial f}{\partial u_1}(a) \right| \geq 1.
$
The conditions of Theorem \ref{TFIU} are satisfied.
As a consequence, we may state the
\begin{prop}\label{TFIU-f}
There exist real numbers  $r\in (0,d),\rho\in (0,r)$ such that , for every  
$a =(a_1,\check{a}) \in \Omega_{\beta/3}$, there exist  a neigbourhood
$U= U_a= B(a_1,r)\times B(\check{a},r)\subset \Omega_{2\beta/3}$ of $a$, an open
 set  $V= V_a$ containing  $\check{a}$ and such that
 $\check{a}\in  B(\check{a},\rho)\subset V\subset B(\check{a},r)$ and 
a function $g_a : V\rightarrow \R$,   $C^2$,  satisfying
\begin{equation}\label{TFI-appli}
\big( (u_1,\check{u}) \in U_a,\ f(u_1,\check{u}) =f(a) \big )
\Longleftrightarrow
\big( \check{u}\in V_a,\ u_1= g_a(\check{u}) \big).
\end{equation}
The differential of  $g_a$ on $V_a$ is given by:
\begin{equation}\label{diffg}
Dg_a(\check{u})= - ( \partial_1f( g_a(\check{u}),\check{u}))^{-1}
( \partial_{\check{\ }} f( g_a(\check{u}),\check{u})),
\end{equation}
where we denote by  $\partial_1$ (resp.  $  \partial_{\check{\ }}$)  the
 partial derivative with respect to the first variable (resp.  the derivative
 with respect to the other variables).  
\end{prop}


\subsection{Uniform enlargement of the open sets}
  
We begin by proving that the auxiliary transformation $S$, defined by 
\eqref{Esse}, is diffeomorphic. This implies that the $T_j^{\theta}$ 
are  diffeomorphic too on convenient domains. This property 
is often useful in what follows. 
We then define the open sets $V_j^{\theta}$, which are the original
open sets $U_j^{\theta}$ slightly and uniformly enlarged, this uniformity being
a consequence of the 
uniform implicit Theorem \ref{TFIU}.
The remaining part of the paragraph is devoted to 
 the definition of the real number $r_0$, which gives condition (PE1). 

\begin{lem}\label{S-diff}
The $C^2$ application $S$  is diffeomorphic  on $\Omega_{\beta}$ onto its image.
\end{lem}

\noindent{\bf Proof.}
Take  $u,v\in \Omega_{\beta}$ such that $S(u)=S(v)$. 
The $k-1$ first coordinates of $S$ ensure that $u_2=v_2, \dots u_k=v_k$. 
Moreover,
the last coordinates of $S(u)$ and $S(v)$ give
$$
\Phi_0(u_1, \gamma^{-1} u_2,\dots, \gamma^{1-k}u_k)=
\Phi_0(v_1, \gamma^{-1} u_2,\dots, \gamma^{1-k}u_k).
$$
Define $g : t\mapsto \Phi_0(t, \gamma^{-1} u_2,\dots, \gamma^{1-k}u_k)$ on
 $(-L-\beta, L+\beta)$. Thanks to the hypotheses \eqref{cond}  satisfied by
$\Phi_0$, $g$ is strictly monotonous. Hence the equality only takes place when 
  $u_1=v_1$ and  $\Phi_0$ is  injective on  $ \Omega_{\beta}$.\\

The differential of $S$ is
\begin{equation}\label{diffS}
DS(u)=
\left(
\begin{array}{ccccccccc}
0& \gamma^{-1} & 0 & \dots & 0\\
0& 0&\gamma^{-1} & \dots & 0\\
0 & 0 & \dots&
0& \gamma^{-1} \\
\displaystyle
\gamma^{k-1}\frac{\partial \Phi_0}{\partial u_1}(\Gamma (u)) &
\displaystyle
\gamma^{k-2}\frac{\partial \Phi_0}{\partial u_2}(\Gamma (u)) & \dots 
\end{array}
\right)
\end{equation}
Its determinant is never equal to zero, which, along with the injectivity,
 proves Lemma \ref{S-diff}.
\hfill $\square$\\

{\bf Definition of the open sets  $V_j^{\theta}$} \\ 
Recall that $d$ is defined by  \eqref{def_d} and that the smaller radius 
$\rho$ given by Proposition \ref{TFIU-f} is $<d$. One  imposes 
\begin{equation}\label{epsm}
\eps_m=\frac{\rho}{3}
\end{equation}
One then sets,  if $U_j^{\theta}\neq \emptyset$, 
\begin{equation}\label{Vjomega}
V_j^{\theta}= B_{\eps_m}(U_j^{\theta})
\end{equation}
and one extends the definition of the application of  $T_{j\theta}$ to
$V_j^{\theta}$. According to Lemma
\ref{S-diff}, $T_{j}^{\theta}$ is diffeomorphic on 
$V_j^{\theta}$ onto its image. \\

We now treat the condition  (PE1) of \cite{SAU}.

\begin{lem}\label{PE1}
There exists $r_0$  such that, for every $\theta \in \R$, for every index $j$ 
 such that $U_j^{\theta}\neq \emptyset$, 
$$
B_{r_0}(S(U_j^{\theta})) \ \subset \ S(V_j^{\theta}).
$$
\end{lem}
\noindent{\bf Proof.}
We begin by proving the existence of a real number  $r_0>0$  such that, 
\begin{equation}\label{PE1-ponct}
\forall x\in \overline{\Omega_{2\beta/3}},\quad B(S(x),r_0)
 \subset S(B(x,\eps_m)).
\end{equation}
According to Lemma \ref{S-diff}, $S$ is diffeomorphic and hence homeomorphic 
on  $\Omega_{\beta}$.  Therefore, for every $x \in  \overline{\Omega_{2\beta/3}}$,
  $S(B(x,\eps_m))$ is an open set containing   $S(x)$. Denote by
\begin{equation}\label{Rdex}
R(x)=\sup\{ r>0 \ :\  B(S(x),r) \ \subset \ S(B(x,\eps_m)) \}.
\end{equation}
Suppose that $\inf_{  \overline{\Omega_{2\beta/3}}}R(x)=0$. There exists a sequence
$(x_n)_{n\in \N}\in \overline{\Omega_{2\beta/3}} ^{\N}$ such that $(R(x_n))_n$
converges to $0$. Since $\overline{\Omega_{2\beta/3}}$ is a compact set, we
can assume that  $(x_n)_{n\in \N}$ converges. Let 
 $x_{\infty} \in \overline{\Omega_{2\beta/3}} $ be its limit point. Let 
\begin{equation*}
R'=\sup \left\{ r>0 \ :\  B(S(x_{\infty}),r) \ \subset \ S \left( B\left(x_{\infty},
\frac{\eps_m}{2}\right) \right) \right\}.
\end{equation*}
The function  $S$ is continuous at  $x_{\infty}$. Let  $\eps>0$ be such that:
$\displaystyle
\eps < \frac{R'}{3}.
$
There exists  $\eta>0$ such that 
$\displaystyle
\Big(\Vert y-x_{\infty}\Vert < \eta \ {\rm and} \
 y \in \overline{\Omega_{2\beta/3}}
\Big) \ \Longrightarrow \ \Vert S(y)-S(x_{\infty})\Vert < \eps.
$
Now take 
 $ y\in   \overline{\Omega_{2\beta/3}} $ such that
$\displaystyle \Vert y-x_{\infty}\Vert < \min\left(\eta, \frac{\eps_m}{2}\right).
$
Then, using the triangle inequality, one sees that
$$
B\left(S(y), \frac{R'}{3}\right) \subset B\left( S(x_{\infty}), 2\frac{R'}{3}  \right).
$$
By definition of $R'$, we obtain
$$
B \left(S(y), \frac{R'}{3}\right) \subset B\left( S(x_{\infty}), 2\frac{R'}{3}  \right)\subset 
S\left(B\left(x_{\infty}, \frac{\eps_m}{2}\right)\right).
$$
This yields
$$
S \left(B\left(x_{\infty}, \frac{\eps_m}{2}\right)\right)\subset  S\left(B\left(y,\eps_m\right)\right)
$$
since $B(x_{\infty}, \frac{\eps_m}{2})\subset  B(y,\eps_m)$ 
(by the triangle inequality and remembering that $\Vert x_{\infty}- y\Vert 
<\frac{\eps_m}{2}$).  
To sum up, for $y\in  \overline{\Omega_{2\beta/3}}$ satisfying 
$\Vert y-x_{\infty}\Vert < \min(\eta, \frac{\eps_m}{2}),$
$$
B \left(S \left(y\right), \frac{R'}{3}\right) \subset S \left(B(y,\eps_m) \right). 
$$
This holds in particular  for  $x_n$, provided $n$ is large enough. Then 
$R(x_n) \geq \frac{R'}{3}$ and does not vanish.  This proves 
 \eqref{PE1-ponct}.

We now may prove the lemma.
Let $y\in B_{r_0}(S(U_j^{\theta}))$. There exists  $x\in U_j^{\theta}$ such that 
$y\in B(S(x), r_0)$. Hence,  by \eqref{PE1-ponct}, 
$y\in  S(B(x,\eps_m))\subset S(V_j^{\theta})$.
Finally
$$
B_{r_0}\left(S\left(U_j^{\theta}\right)\right) \subset S \left(V_j^{\theta}\right).
$$
This concludes the proof of Lemma \ref{PE1}.\hfill $\square$\\


\subsection{Lipschitz continuity of the Jacobian determinant}

We go on with the regularity properties of the transformations $T_j^{\theta}$, 
which come from the properties of the auxiliary function $S$. We prove that the 
Jacobian determinant satisfies the following  Lipschitz property:

\begin{prop}\label{hoelder}
There exist $K > 0$ and $\tilde{\eps} > 0$ such that, for all 
$\theta \in \R$, $j \in \Z$
 satisfying $U_j^\theta \ne \emptyset$, for all $\varepsilon \in (0,
 \tilde{\eps})$, $z \in S(\Omega_{\beta/3})$ and
 $x, y \in B_\varepsilon (z)$,
\begin{equation}\label{detholder}
\abs{\left(\det DS^{-1}(x)\right) - \left(\det DS^{-1}(y)\right)} \le
 \frac K2 \abs{\left(\det DS^{-1}(z)\right)} \norm{x - y} .
\end{equation}
\end{prop}

{\bf Proof.}
For $z \in S(\Omega_{\beta/3})$ define
\[
\begin{matrix}
	g_z : & S(\Omega_{\beta}) & \longrightarrow & \R \\
	& x & \longmapsto & \left(\det DS^{-1}(z)\right)^{-1} \left(\det DS^{-1}(x)\right)
\end{matrix}
\]

From  \eqref{diffS} we get that $\det DS(S^{-1}(x)) = 
(-1)^{k-1} \dfrac{\partial\Phi_0}{\partial x_1} (\Gamma S^{-1}(x))$.
 Hence for all  $x \in S(\Omega_{\beta})$, 
 $g_z (x) = \left(\dfrac{\partial\Phi_0}{\partial x_1} (\Gamma S^{-1}(z))\right)
 \left(\dfrac{\partial\Phi_0}{\partial x_1} (\Gamma S^{-1}(x))\right)^{-1}$ and
\begin{align}
Dg_z(x) &= - \left(\dfrac{\partial\Phi_0}{\partial x_1} 
(\Gamma S^{-1}(z))\right) \left(\dfrac{\partial\Phi_0}{\partial x_1} 
(\Gamma S^{-1}(x))\right)^{-2} D \left(x \mapsto \dfrac{\partial\Phi_0}
{\partial x_1} (\Gamma S^{-1}(x))\right). \notag 
\end{align}
Using conditions \eqref{cond} and the  $C^2$ regularity of $\Phi_0$ gives that,
 for $x \in S(\overline{\Omega_{2\beta/3}})$,
\begin{align*}
\norm{D g_z (x)}_{op} &\le \left|\dfrac{\partial\Phi_0}{\partial x_1} (\Gamma S^{-1}(z))\right| \left|\dfrac{\partial\Phi_0}{\partial x_1} (\Gamma S^{-1}(x))\right|^{-2} \norm{D \dfrac{\partial\Phi_0}{\partial x_1} (\Gamma S^{-1}(x))}_{op} \norm{\Gamma}_{op} \norm{DS^{-1} (x)}_{op} \\
&\le \sup \left|\dfrac{\partial\Phi_0}{\partial x_1}\right| C_1^{-k + 1} C_2^{-1} \sigma^{-k} \sup \norm{\nabla \dfrac{\partial\Phi_0}{\partial x_1}} \gamma^{-k+1} \sup \norm{DS^{-1}}_{op} \eqqcolon K/2 
\end{align*}
where the upper bounds are taken on  $\left[-L - \frac{2\beta}3,
 L + \frac{2\beta}3\right]^k$, except for the last one, taken on 
 $S(\overline{\Omega_{2\beta/3}})$. Here $\norm{\cdot}_{op}$ denotes the 
operator norm associated with the Euclidean norm.\\

Choose $\tilde{\eps} =  d\left( S\!\left( \overline{\Omega_{\beta/3}} \right)\!, 
\complement S(\Omega_{2\beta/3}) \right)$ (this is $>0$). If
 $\eps <\tilde{\eps}$ and  $x, y \in B_\eps (z)$,
 the segment  $[x, y]$ is included into 
$B_\eps (z) \subset S \left( \Omega_{2\beta/3} \right)$.
Then, by the fundamental theorem of differential calculus, 
\[
\abs{g_z (x) - g_z (y)} \le \frac K2 \norm{x - y}.
\]
Dividing both terms by  $\abs{\left(\det DS^{-1} (z)\right)}$ yields 
 \eqref{detholder}.
\qed

Since, up to an additive constant,  the applications  $T_{j}^{\theta}$ coincide
 with $S$ on the open sets  $V_j^\theta$, we obtain the following consequence:
\begin{cor}
There exist  $K > 0$ and $\tilde{\eps} > 0$ such that, for all 
 $\theta \in \R$, $j \in \Z$ satisfying $U_j^\theta \ne \emptyset$,
 $\varepsilon \in (0, \tilde{\eps})$, $z \in
 T_{j}^{\theta}(V_j^\theta)$ and $x, y \in B_\varepsilon (z) \cap T_{j}^{\theta} (V_j^\theta)$, 
\[
\abs{\det D(T_{j}^{\theta})^{-1} (x) - \det (DT_{j}^{\theta})^{-1} (y)} \le K
 \abs{\det (DT_{j}^{\theta})^{-1} (z)} \varepsilon
\]
\end{cor}

The condition (PE2) in \cite{SAU} holds. 

\subsection{Geometry of the boundaries}

This paragraph concludes the series of proofs concerning the conditions 
listed in Theorem \ref{unif}.\\
Instead of studying the expression which appears in \eqref{eta-G}, 
we will rely on Lemma 2.1 of \cite{SAU}, which gives a simpler 
sufficient condition. Both conditions are of geometric type: the first 
one is concerned
with volumes of balls centered on the boundary, the second one, with the
 maximal number of 
crossings at one point (giving the number of open sets that the said balls
 intersect). The condition stated in \eqref{eta-G} is satisfied 
if the limit of a quantity (that  we will define and  study below) is $<1$. 
Then, for a sufficiently small parameter, the quantity itself is $<1$. 

We will follow
the proof of \cite {SAU} and check that the ``sufficiently small'' character 
does not depend too severely on the perturbation $\theta$.
We need to give an upper bound for the volumes which appear in \eqref{eta-G}.
To that end, we rectify the boundaries of the open sets and are brought back
to the volume of the portion of a ball comprized between two parallel 
hyperplanes. We rely on Theorem \ref{TFIU-f},  on the uniform open charts 
defined  below in \eqref{carte} and which we  use, too,  to prove that the 
boundaries are null sets. 

\noindent{\it Local charts, negligeability of the boundaries}\\ 
We are concerned with the boundaries of the open sets  $U_j^{\theta}$ defined by
\eqref{ouv-front}. There are portions of hyperplanes, which we do not need to
rectify and the  $\Sigma^{\theta}_j $. Recall that the function $f$ is defined by
$f=\Phi_0\circ \Gamma$ (see \eqref{def-phi-om},\eqref{Gamma-dilate},
\eqref{f-TFI}). The distance $d$ is defined in \eqref{def_d}, the real numbers 
$0<\rho<r$ are given by Proposition  \ref{TFIU-f}, of which we will use the 
notations  $U_a, V_a$ and $ g_a$.\\

\noindent{\it  Definition of the charts: }\\
For $a=(a_1,\check{a})\in \Omega_{\beta/3}$, set
\begin{equation}\label{Uprima}
U'_a:=B(a_1,\rho) \times B(\check{a},\rho).
\end{equation}
We define the chart by
\begin{equation}\label{carte}
\forall u=(u_1,\check{u})\in U'_a,\quad 
\Phi_a(u)=  (u_1-g_a(\check{u}), \check{u}).
\end{equation}
According to Proposition \ref{TFIU-f}, $B(\check{a},\rho)\subset V_a \subset 
B(\check{a},r)$ and $g_a$, defined  on $V_a$, is defined on the smaller
set $B(\check{a},\rho)$.

The application $\Phi_a $ is $C^2$ on $U'_a$ and clearly injective. Moreover, 
\begin{equation}\label{jac}
D\Phi_a(u)= 
\left( 
\begin{array}{lllll}
 1& -\partial_2g_a(\check{u}) &   -\partial_3 g_a(\check{u})  & \dots &  -\partial_k g_a(\check{u})   \\
0 & 1 & 0 & \dots & 0 \\
0 &0&1&\dots  \\
\dots \\
0 & 0 & \dots & 0 & 1
\\\\
\end{array}
\right).
\end{equation}
The inverse function theorem then proves that 
 $\Phi_a$ is diffeomorphic on  $U'_a$ onto its image. 

Denote by  $H$ the hyperplane  $H=\{ (0,\check{v}),\check{v} \in \R^{k-1}\}$.
Let  $u\in U'_a$ satisfy  $f(u)=f(a)$.  Then, according to  \eqref{TFI-appli},
 $\check{u}\in V_a$  and 
$u_1= g_a(\check{u})$. Therefore $\Phi_a(u)=(0,\check{u})\in H$. To sum up: 
\begin{equation}\label{rectif}
\Phi_a(\{ u\in U'_a : f(u)=f(a)\})  \ \subset \ H.
\end{equation}
This boundary is locally rectified.

 The compact set  $\Sigma_j^\theta$  is covered by a finite
number of charts. Since each local chart  $\Phi_a$ has determinant $1$, 
the volumes are conserved. 
This gives the following result:
 
\begin{lem}\label{front-neg}
The sets  $\Sigma_j^\theta \ (j\in \Z,\theta \in \R)$   are $C^2$ 
 and have Lebesgue measure zero.
\end{lem}


\noindent{\it Number of crossings}\\
Consider a given point $x\in \R^k$ and a small neigbourhood $\cal O$ of $x$.
We wish to give an upper bound for the number of boundaries of the open sets
$U_j^{\theta}$ which intersect $\cal O$. Some boundary parts are the sides of the
parallelepiped $\Omega$, of which  at most $k$ intersect  $\cal O$, if 
one is near a corner.
We prove here that, if $\cal O$ is sufficiently small, it intersects at most 
 one slanted boundary $\Sigma_j^k$ (whatever the value of $\theta$).  
This slanted boundary counts as two, because $\cal O$ intersects two open
sets of the partition and that makes two terms in the sum of volumes 
we must consider (see \eqref{eta-G} and \eqref{eta-GG}). The maximal number
 of crossings $Y$ is then $k+2$. Except for a countable number of parameters
 $\theta$, this number is even $k+1$. 

\begin{lem}\label{dist-front}
Let $\theta\in \R$, let $j,j'$ be two different integers such that 
$\Sigma_{j}^{\theta}$ and  $\Sigma_{j'}^{\theta}$ are not empty. Take one point in 
each of them: 
$x \in \Sigma_{j}^{\theta}$, $x'\in \Sigma_{j'}^{\theta}$.\\
Let $M^1_{\Phi_0}$ and $M_{\Gamma}$  be positive real numbers such that
$$
\forall u\in [-L,L]^k,\ \Vert D\Phi_0(u)\Vert \leq M^1_{\Phi_0} \ {\rm and}\
\Vert \Gamma\Vert \leq M_{\Gamma}
$$
(this should hold for operator-norms associated with all the norms 
on $\R^k$ that we need to use).Then
$$
\Vert x-x'\Vert \geq \frac{2L}{M^1_{\Phi_0}\ M_{\Gamma} }.
$$
As a consequence, if
\begin{equation}\label{encore}
\eps_0 < \frac{L}{2  M^1_{\Phi_0}\ M_{\Gamma} },
\end{equation}
the ball $B_{\eps_0}(x)$ can't intersect both bands
$B_{ \eps_0}(\Sigma_j^{\theta})$  and $B_{ \eps_0}(\Sigma_{j'}^{\theta}).$
\end{lem}

\noindent{\bf Proof.}\\
The hypotheses imply that
$$
\Phi_0(\Gamma(x))= -\theta + (2j-1)L,\quad \Phi_0(\Gamma(x'))= -\theta + (2j'-1)L.
$$
Applying the fundamental theorem of differential calculus gives 
$$
|\Phi_0(\Gamma(x))- \Phi_0(\Gamma(x'))|\leq 
\sup_{y\in[x,x']}\Vert D\Phi_0(y)\Vert \Vert \Gamma\Vert \Vert x-x'\Vert.
$$
Hence
$$
|\Phi_0(\Gamma(x))- \Phi_0(\Gamma(x'))|
= | 2(j-j')L | \leq  M^1_{\Phi_0}\ M_{\Gamma} 
\Vert x-x'\Vert.
$$
But $|j-j'|$ is greater than $1$, which gives the first result. The last point
 is then a consequence of the triangle inequality.  \hfill $\square$\\


\noindent{\it A uniform version of Lemma  2.1 of \cite{SAU}}\\
Let  $\theta\in \R$ be a fixed perturbation, let  $\ph_{\theta}$ and $T_{\theta}$,
be associated with $\theta$ by  \eqref{Tomega} and \eqref{T-et-S}.
We defined  $\eps_m>0$ and the neigbourhoods $V_j^{\theta}$ of the $U_j^{\theta}$ 
by \eqref{epsm} and \eqref{Vjomega}:
\begin{equation*}
\eps_m=\frac{\rho}{3}, \quad V_j^{\theta}= B_{\eps_m}(U_j^{\theta})\subset 
\Omega_{\beta/3}.
\end{equation*}
Lemma \ref{PE1} gives a positive real  number $r_0$  such that, 
 for all  $\theta' \in \R$ and all $j$ such that  
 $U_j^{\theta'}\neq \emptyset$, 
$$
B_{r_0} (S(U_j^{\theta'})) \ \subset \ S(V_j^{\theta'}).
$$

Let $x\in \R^k$, let $\eps,\eps_0$ satisfy
\begin{equation}\label{epszero}
0<\eps<\eps_0<\min \left(\eps_m,r_0,\frac{L}{2  M^1_{\Phi_0}\ M_{\Gamma} },\tilde{\eps} \right)
\end{equation}
with the notations of Lemma \ref{dist-front} and Proposition \ref{hoelder}.
 Recall that $T_{j}^{\theta}$ is 
considered as defined on  $V_j^{\theta}$, although its expression is well defined 
on the much broader set $\Omega_{\beta}$.\\

To check the condition (PE5), we must find an upper bound for the volume
of the following set. 
\begin{equation}\label{horreur}
(T_{j}^{\theta})^{-1} (B_{\eps}(\partial T_{j}^{\theta}(U_j^{\theta})))\ 
 \cap B_{(1-\frac{1}{\sqrt{\sigma}} )\eps_0}(x).
\end{equation}

According to Lemma \ref{S-diff}, $T_{\theta}$ is diffeomorphic on  $V_j^{\theta}$ 
onto its image. Moreover, it is expanding by Proposition \ref{dilatance}.
As in \cite{SAU}, we may state the following inclusions:
\begin{equation}\label{incl-1}
(T_{j}^{\theta})^{-1} (B_{\eps}(\partial T_{j}^{\theta}(U_j^{\theta})))\subset 
B_{\frac{1}{\sqrt{\sigma}} \eps}(\partial U_j^{\theta})\subset V_j^{\theta} \subset
 B_{\frac{1}{\sqrt{\sigma}} \frac{d}{3}}(\Omega).
\end{equation}
Formula \eqref{def_d} which defines $d$, implies that the set
\eqref{horreur} is empty if  $x\notin \Omega_{\beta/3}$.

The boundary of an open set $ U_j^{\theta}$ consists of plane sets (sides of a
 parallelepiped) and of slanted boundaries like  $\Sigma_j^{\theta}$.
If we denote by $F_s$ these different kinds of sets, we obtain
$$
 (T^{\theta})^{-1} (B_{\eps}(\partial T_{\theta}(U_j^{\theta})))\ 
 \cap B_{(1-\frac{1}{\sqrt{\sigma}} )\eps}(x) \subset  \bigcup
 B_{\frac{1}{\sqrt{\sigma}} \eps}(F_s)\cap  B_{(1-\frac{1}{\sqrt{\sigma}} )\eps_0}(x) ,
$$
and we will treat below each case (plane or not) separately. \\

{\bf First case: $F_s$ is included in a hyperplane} \\
The set $B_{\frac{1}{\sqrt{\sigma}} \eps}(F_s)\cap B_{(1-\frac{1}{\sqrt{\sigma}})\eps_0}(x)$
has maximal volume when $F_s$ contains the center $x$ of the ball. In this case,
the volume is smaller than that of a cylinder based on a $k-1$ dimensional ball,
which gives:
\begin{equation}\label{vol-simple}
m\left( B_{\frac{1}{\sqrt{\sigma}} \eps}(F_s)\cap  B_{(1-\frac{1}{\sqrt{\sigma}} )\eps_0}(x)
 \right)
\leq \underbrace{2\frac{\eps}{\sqrt{\sigma}}}_{ {\rm height}} \times 
\underbrace{\gamma_{k-1}
\left((1-\frac{1}{\sqrt{\sigma}})\eps_0\right)^{k-1}}_{{\rm volume \ of \ the \ 
ball}}.
\end{equation}

{\bf Second case: $F_s= \Sigma_j^{\theta}$}\\
We compute an upper bound for the volume of the set below
(supposed to be non empty):
\begin{equation}\label{inter}
B_{\frac{1}{\sqrt{\sigma}} \eps}(\Sigma_j^{\theta})\cap 
 B_{(1-\frac{1}{\sqrt{\sigma}} )\eps_0}(x) .
\end{equation}

We will use  some auxiliary points and we need to locate them exactly.

\noindent{\it Location  of auxiliary points $u,a,y,z$}
\\
Let $u$ belong to the intersection. It is in the first set, which gives
 an  $a \in \Sigma_j^{\theta}$ (hence in $\Omega\subset \Omega _{\beta/3}$ ) such
 that $\Vert u-a\Vert < \frac{1}{\sqrt{\sigma}} \eps$.
Recall that, thanks to \eqref{rectif}, $\Phi_a(\Sigma_j^{\theta}\cap U'_a)\subset H$.
\\
It is in the ball too, hence 
$\Vert u-x\Vert < (1-\frac{1}{\sqrt{\sigma}} )\eps_0$.

Therefore
$$
\Vert a-x\Vert < \frac{1}{\sqrt{\sigma}}
 \eps+(1-\frac{1}{\sqrt{\sigma}} )\eps_0 \leq \eps_0
<\eps_m=\frac{\rho}{3}< \frac{d}{3}.
$$
Then $x\in \Omega_{\beta/3}$ (thanks to $d$) and $x$ is in the chart open set 
$U'_a$ (see \eqref{Uprima}) associated with  $a\in \Sigma_j^{\theta}$
 (thanks to $\rho$).\\
 If $y$ is in $B_{(1-\frac{1}{\sqrt{\sigma}} )\eps_0}(x) $, 
$\Vert a-y\Vert < 2 \eps_0.$  Hence 
$$
B_{(1-\frac{1}{\sqrt{\sigma}} )\eps_0}(x)\subset B(a, 2 \eps_0).
$$
This is a subset of $U'_a$ and of $ \Omega_{\beta/3}$ because $2\eps_0 < 2 
\rho/3 < 2 d/3$.

And if $y\in  B_{\frac{1}{\sqrt{\sigma}} \eps}(\Sigma_j^{\theta})\cap 
 B_{(1-\frac{1}{\sqrt{\sigma}} )\eps_0}(x) $,
there exists  $z\in \Sigma_j^{\theta}$ such  that 
$\Vert y-z \Vert <\frac{1}{\sqrt{\sigma}} \eps$.
Then $\Vert a-z\Vert < 3\eps_0$  and  $z$ is in $U'_a$ too.
We will use these remarks in Lemma \ref{TaylorRI}.  \\

Remark that,  if $O$ is an open subset of $U'_a$, 
$ \displaystyle m(\Phi_a(O))= m(O).$
This  is because  $\Phi_a$ preserves the volume, as was already used to obtain
Lemma \ref{front-neg}. For the same reason, the linear transformation 
$D\Phi_a(x)$  preserves the volumes too :
\begin{equation}\label{aidebis}
m(D\Phi_a(x)^{-1}\Phi_a(O))= m(O).
\end{equation}

We are brought back to the computation of the following measure
$$
m\left( D\Phi_a(x)^{-1}\Phi_a(B_{\frac{1}{\sqrt{\sigma}} \eps}
(\Sigma_j^{\theta})\cap U'_a)\cap
 D\Phi_a(x)^{-1}
 \Phi_a( B_{(1-\frac{1}{\sqrt{\sigma}} )\eps_0}(x) )\right).
$$
To this aim we  have to locate the open set relatively to a convenient
 hyperplane.

We first establish some technical results. 
\begin{lem}\label{TaylorRI}
Let the points $x,a$ be as above and take $u,v$ in $U'_a$. Then
we can write
\begin{equation}\label{aide-form}
\begin{array}{llll}
\displaystyle
 \Vert D\Phi_a(x)^{-1} \Phi_a(u) - D\Phi_a(x)^{-1} \Phi_a(x) \Vert
\leq \Vert u-x\Vert + M \Vert u-x\Vert^2.\\ \\
\displaystyle
\Vert  D\Phi_a(x)^{-1}\Phi_a(u)-  D\Phi_a(x)^{-1}\Phi_a(v) \Vert \leq
 \Vert u-v\Vert + M \Vert v-x \Vert   \Vert u-v\Vert + M
\Vert u-v\Vert ^2 
\end{array}
\end{equation}
with a constant  $M>0$, depending only on the bounds of $D\Phi_0, D^2\Phi_0$ 
and $(\partial_1\Phi_0)^{-1}$ on $\overline{\Omega_{2\beta/3}}$.
\end{lem}
\noindent{\bf Proof.}
We apply the Taylor integral formula:
$$
\begin{array}{lll}
\Phi_a(u)-\Phi_a(x) &\displaystyle = D\Phi_a(x)\cdot (u-x)+ 
\int_0^1 (1-t)D^2\Phi_a(x+t(u-x))(u-x)^2 \ dt \\\\
&\displaystyle = D\Phi_a(x)\cdot \Big[ (u-x)+ 
\int_0^1 (1-t) D\Phi_a(x)^{-1}D^2\Phi_a(x+t(u-x))(u-x)^2 \ dt\Big]. \\
\end{array}
$$
Therefore,
\begin{equation}\label{a}
 D\Phi_a(x)^{-1} \Phi_a(u) - D\Phi_a(x)^{-1} \Phi_a(x) =
 (u-x)+ 
\int_0^1 (1-t) D\Phi_a(x)^{-1}D^2\Phi_a(x+t(u-x))(u-x)^2 \ dt ,
\end{equation}
which will give the first inequality. To get the second one, we write
$$
\Phi_a(u)-\Phi_a(v) = 
  D\Phi_a(v)\cdot (u-v)+ 
\int_0^1 (1-t)D^2\Phi_a(v+t(u-v))(u-v)^2 \ dt 
$$
and develop the term $ D\Phi_a(v)$. This gives
$$
\Phi_a(u)-\Phi_a(v) = 
  D\Phi_a(x)\cdot (u-v)+  \int_0^1 D^2\Phi_a(x+s(v-x))(v-x) \ ds\cdot (u-v)+
\int_0^1 (1-t)D^2\Phi_a(v+t(u-v))(u-v)^2 \ dt .
$$
Compose by $  D\Phi_a(x)^{-1}$: 
\begin{equation}\label{b}
\begin{array}{lll}
\displaystyle
 D\Phi_a(x)^{-1}\Phi_a(u)
-  D\Phi_a(x)^{-1}\Phi_a(v) & =  (u-v)  \\  \\
&\displaystyle+
\int_0^1 D\Phi_a(x)^{-1} D^2\Phi_a(x+s(v-x))(v-x) \ ds\cdot (u-v) \\ \\
&\displaystyle+
\int_0^1 (1-t) D\Phi_a(x)^{-1}  D^2\Phi_a(v+t(u-v))(u-v)^2 \ dt .
\end{array}
\end{equation}
All points appearing in the integrals belong to  $U'_a $, which is a 
convex subset of $\Omega_{2\beta/3}$. The first, second differentials and the
 inverses appearing in the integrals are bounded, in norm, on
 $\Omega_{2\beta/3}$. If we call $M$ 
a common upper bound, we obtain, from \eqref{a} and \eqref{b} respectively, the
first and the second inequality of \eqref{aide-form}.
This completes the proof of Lemma \ref{TaylorRI}.  \hfill $\square$\\

We now turn to the set $E$ defined by
\begin{equation}\label{E2}
 E= D\Phi_a(x)^{-1} \Phi_a \left( B_{(1-\frac{1}{\sqrt{\sigma}} )\eps_0}(x) \right).
\end{equation}

According to Lemma \ref{TaylorRI}, if $y\in B_{(1-\frac{1}{\sqrt{\sigma}})\eps_0}(x)$, 
\begin{equation*}
 \Vert D\Phi_a(x)^{-1} \Phi_a(y) - D\Phi_a(x)^{-1} \Phi_a(x) \Vert
\leq \Vert y-x\Vert(1 + M \Vert y-x\Vert).
\end{equation*}

Then
\begin{equation}\label{E2inclus}
D\Phi_a(x)^{-1} \Phi_a(  B_{(1-\frac{1}{\sqrt{\sigma}} )\eps_0}(x) )
\subset  B( D\Phi_a(x)^{-1} \Phi_a(x) , R) := F,
\end{equation}
with 
$$R={\left(1-\frac{1}{\sqrt{\sigma}}\right)\eps_0}
\left(1+M{\left(1-\frac{1}{\sqrt{\sigma}}\right)\eps_0}\right),$$
since
$\Vert y-x\Vert \leq (1-\frac{1}{\sqrt{\sigma}} )\eps_0$.\\

Now suppose that $y$ is in  
$ B_{\eps\frac{1}{\sqrt{\sigma}}  }(\Sigma_j^{\theta})\cap
  B_{(1-\frac{1}{\sqrt{\sigma}} )\eps_0}(x)$. 
We want to get an upper bound for the distance between the point
 $D\Phi_a(x)^{-1}\Phi_a(y)$ and the set  $D\Phi_a(x)^{-1} H$, where $H$
is the hyperplane already defined by
 $H=\{ (0,\check{v}),\check{v} \in \R^{k-1}\}$.
 Let 
$z\in \Sigma_j^{\theta}$ be such that $\Vert y-z\Vert 
<\eps\frac{1}{\sqrt{\sigma}} $. Then $z \in U'_a$ and by Lemma 
\ref{TaylorRI}, 
$$
\Vert D\Phi_a(x)^{-1}\Phi_a(y)-D\Phi_a(x)^{-1}\Phi_a(z)\Vert \leq 
\Vert y-z\Vert  (1 + M \Vert  y-z \Vert  +M \Vert x-z\Vert).
$$
Remark that
$$
\Vert x-z\Vert\leq \Vert x-y\Vert + \Vert y-z\Vert
\leq{\left(1-\frac{1}{\sqrt{\sigma}} \right)\eps_0} + \eps\frac{1}{\sqrt{\sigma}},
$$
with $\eps \leq\eps_0$. Hence 
$$
d(D\Phi_a(x)^{-1}\Phi_a(y),D\Phi_a(x)^{-1} H ) \leq 
\eps\frac{1}{\sqrt{\sigma}}  \left(1+ M \eps_0\left(1+\frac{1}{\sqrt{\sigma}}\right) \right).
$$

The set  $H' := D\Phi_a(x)^{-1} H$ is a hyperplane too. Therefore
\begin{equation}\label{E1inclus}
D\Phi_a(x)^{-1}  \Phi_a(B_{\frac{1}{\sqrt{\sigma}} \eps}(\Sigma_j^{\theta})\cap 
 B_{(1-\frac{1}{\sqrt{\sigma}} )\eps_0}(x) )\ \subset \
B\left(D\Phi_a(x)^{-1} H, \eps\frac{1}{\sqrt{\sigma}} 
 (1+ M \eps_0(1+\sigma^{-1/2}) ) \right) :=F_1.
\end{equation}

The  inclusions \eqref{E1inclus} and \eqref{E2inclus} give:
$$
D\Phi_a(x)^{-1}\Phi_a (B_{\frac{1}{\sqrt{\sigma}} \eps}(\Sigma_j^{\theta})\cap 
B_{(1-\frac{1}{\sqrt{\sigma}} )\eps_0}(x))\subset  F\cap F_1
$$
We have seen, too, that this set
$D\Phi_a(x)^{-1}\Phi_a (B_{\frac{1}{\sqrt{\sigma}} \eps}(\Sigma_j^{\theta})\cap
  B_{(1-\frac{1}{\sqrt{\sigma}} )\eps_0}(x))$ has the same volume as 
$$
B_{\eps \frac{1}{\sqrt{\sigma}} }(\Sigma_j^{\theta})
 \cap B_{(1-\frac{1}{\sqrt{\sigma}} )\eps_0}(x),
$$
which is the set we are interested in. 
It implies
\begin{equation}\label{mes}
m\left( B_{\eps \frac{1}{\sqrt{\sigma}} }(\Sigma_j^{\theta}) 
 \cap B_{(1-\frac{1}{\sqrt{\sigma}} )\eps_0}(x)\right) \leq
m(F_1\cap F).
\end{equation}

Then  $F_1\cap F$ is the intersection of a ball with the space comprised
 between two hyperplanes parallel to $H'$ and situated at the same distance. 
The volume in maximal when the center is on $H'$ and it is less than 
\begin{equation}\label{volume}
\begin{array}{llll}
\displaystyle
 2 \underbrace{ \eps\frac{1}{\sqrt{\sigma}}  (1+ M \eps_0(1+\sigma^{-1/2}) )}_{{\rm height}} 
\times \gamma_{k-1} 
\underbrace{\left( {(1-\frac{1}{\sqrt{\sigma}} )\eps_0}(1+M{(1-\frac{1}{\sqrt{\sigma}} )\eps_0}) \right)^{k-1}}_{{\rm radius}=R}\\ \\
= 
\displaystyle
 2 \gamma_{k-1}  \eps\frac{1}{\sqrt{\sigma}} 
\left( {\left(1-\frac{1}{\sqrt{\sigma}} \right)\eps_0} \right)^{k-1}
 \times  \left(1+ M \eps_0 \left(1+\sigma^{-1/2}\right) \right)
\left( \left(1+M{(1-\frac{1}{\sqrt{\sigma}} )\eps_0}\right) \right)^{k-1}
\end{array}
\end{equation}
This upper bound holds for the simpler first case, when the boundary is
a side of $\Omega$ (see \eqref{vol-simple}).

According to \cite{SAU}, we must find an upper bound for
\begin{equation}\label{Geps}
 G(\eps,\eps_0) = \sup_{x\in \R^n} G(x,\eps,\eps_0)
\end{equation} with 
\begin{equation}\label{Gtout}
 G(x,\eps,\eps_0)= \sum_{i}
\frac{m(T_i^{-1}(
 B_\eps(\partial T^{\theta}U_i^{\theta}))
\cap 
B(x,(1-\frac{1}{\sqrt{\sigma}})\eps_0) )}
{m(B(x,(1-\frac{1}{\sqrt{\sigma}})\eps_0))}
\end{equation}

Precisely,  $\eta(\eps_0)$ is defined by 
\begin{equation}\label{eta-GG}
\eta(\eps_0)= \frac{1}{\sqrt{\sigma}} + 2
\sup_{\eps\leq \eps_0} \frac{ G(\eps,\eps_0)}{\eps} \ \eps_0.
\end{equation}

What precedes establishes that the sum (indexed by $i$) of all the volumes
we have to consider is smaller than
$$
Y\times 
 2 \gamma_{k-1}  \eps\frac{1}{\sqrt{\sigma}} 
\left( {\left(1-\frac{1}{\sqrt{\sigma}} \right)\eps_0} \right)^{k-1}
 \times  \left(1+ M \eps_0 \left(1+\sigma^{-1/2}\right) \right)
\left( \left(1+M{\left(1-\frac{1}{\sqrt{\sigma}} \right)\eps_0}\right) \right)^{k-1},
$$
with $Y= k+1$ if no slanted boundary  $\Sigma_j^{\theta}$ meets the corner of 
 $\Omega$, and $Y=k+2$ otherwise. The slanted boundary appears twice in the sum,
for two successive indices $i$. 

With the definitions \eqref{Geps} and \eqref{Gtout} we obtain
$$
 G(x,\eps,\eps_0)\leq
\frac{2 \gamma_{k-1} Y \eps\frac{1}{\sqrt{\sigma}}}{ m\left( B \left( x, \left( 1-\frac{1}{\sqrt{\sigma}} \right) \eps_0 \right) \right) }  \left( {\left(1-\frac{1}{\sqrt{\sigma}} \right)\eps_0} \right)^{k-1}
  (1+ M \eps_0(1+\sigma^{-1/2}) )
\left( \left(1+M{ \left(1-\frac{1}{\sqrt{\sigma}} \right)\eps_0}\right) \right)^{k-1}  ,
$$
which can be written more simply as
$$
 G(x,\eps,\eps_0)\leq
\frac{ 2 \gamma_{k-1} Y}{\gamma_k  (\sqrt{\sigma}-1)\eps_0 } 
 \eps 
 \times  (1+ M \eps_0(1+\sigma^{-1/2}) )
\left( \left(1+M{\left(1-\frac{1}{\sqrt{\sigma}} \right)\eps_0}\right) \right)^{k-1}  .
$$
This does not depend on $x$. Similarly,
$$
\frac{ G(\eps,\eps_0)}{\eps} \ \eps_0\leq 
\frac{ 2 \gamma_{k-1} Y}{\gamma_k  (\sqrt{\sigma}-1) } 
 \times  (1+ M \eps_0(1+\sigma^{-1/2}) )
\left( \left(1+M{\left(1-\frac{1}{\sqrt{\sigma}} \right)\eps_0}\right) \right)^{k-1} ,
$$
which
does not depend on $\eps$. We finally obtain that the function
$\eta$ is smaller than an increasing function $\bar{\eta}$ of $\eps_0$ :
\begin{equation}\label{etabarre}
\begin{array}{lll}\displaystyle
\eta(\eps_0)= \frac{1}{\sqrt{\sigma}} + 2
\sup_{\eps\leq \eps_0}
\frac{ G(\eps,\eps_0)}{\eps} \ \eps_0\\
\displaystyle\leq  \overline{\eta(\eps_0)} :=
\frac{1}{\sqrt{\sigma}} + 
\frac{ 4 \gamma_{k-1} Y}{\gamma_k (\sqrt{\sigma}-1) } 
 \times  (1+ M \eps_0(1+\sigma^{-1/2}) )
\left( \left(1+M{\left(  1-\frac{1}{\sqrt{\sigma}} \right)\eps_0}\right) \right)^{k-1} ,
\end{array}
\end{equation}
The limit of the upper bound $\overline{\eta(\eps_0)}$, as $\eps_0$ goes
 to $0$, is
 $$
\eta_Y(\sigma)=
\frac{1}{\sqrt{\sigma}} + 
\frac{ 4\gamma_{k-1} Y}{\gamma_k  } 
 \frac{1}{\sqrt{\sigma} \ -1},
$$
which is assumed to be $<1$ (see \eqref{eta}).
Therefore, there exists an $\eps _0$ such  that (PE5) is satisfied, 
which proves that all
the conditions of Theorem 5.1 in \cite{SAU}, listed here in Theorem \ref{unif},
 are satisfied. 

\section{Proof of Theorems \ref{av-TPF} and \ref{th-inv-meas} : Properties of the average transfer operator}
\label{sec-avPF}

Seeing as $T_\theta : \Omega \to \Omega$ is nonsingular for every $\theta \in \R$, we can define its (Perron-Frobenius) transfer operator 
$P_\th$.

Our goal is to define an averaged form of the transfer operator:
\begin{equation}\label{PF_moyen}
	\hspace{1em} Pf = \int_\R P_\theta f \d\Pr_{\Theta_0}(\theta) = \int_\R H(\theta) P_\theta f \d\theta,
\end{equation}
for $f \in L^1(\Omega, m)$, where $H$ is a density function. This definition raises the issue of measurability in the integrand, which we shall address before establishing spectral properties for $P$.

\subsection{Measurability and $L^p$-boundedness}

\begin{prop}
Let $f \in L^1 (\Omega, m)$.
\begin{enumerate}
\item The mapping
\begin{equation}\label{forme_PF}
\begin{matrix}
Qf : & \Omega \times \R& \longrightarrow & \R \\
& (y, \theta) & \longmapsto & \displaystyle \sum_{j \in \Z} \1_{T_{j}^{\theta} (U_j^\theta)} (y) \left(f \circ (T_{j}^{ \theta})^{-1}\right) \! (y) \abs{\frac{\partial \Phi_0}{\partial x_1} \left(\Gamma (T_{j}^{ \theta})^{-1} (y)\right)}^{-1}
			\end{matrix}
		\end{equation}
		is measurable.
		\item %
		For all $(y, \theta) \in \Omega \times \R$, $Qf(y, \theta) = P_\theta f(y)$.
		\item The mapping
		\[
		\begin{matrix}
			\bar Q f : & \Omega \times \R& \longrightarrow & \R \\
			& (y, \theta) & \longmapsto & H(\theta) P_\theta f(y)
		\end{matrix}
		\]
		is measurable and integrable over $\Omega \times \R$.
	\end{enumerate}
\end{prop}

\noindent{\bf Proof.}
\begin{enumerate}
	\item For fixed $j \in \Z$,  since the set $\{ (y,\theta)\in \R^k \times \R : y \in T_{j}^{\theta}(U_j^{\theta}) \}$ is open, the map $(y, \theta) \mapsto \1_{T_{j}^{\theta} (U_j^\theta)} (y)$ is measurable. 
 So are the maps $(y, \theta) \mapsto \left(f \circ (T_{j}^{ \theta})^{-1}\right) (y)$ and $(y, \theta) \mapsto \abs{\frac{\partial \Phi_0}{\partial x_1} \left(\Gamma (T_{j}^{ \theta})^{-1} (y)\right)}^{-1}$. $Qf$ is therefore a series of measurable functions. In addition, for fixed $\theta$, only finitely many of the $U_j^\th$ are nonempty, thus for every $(y, \th) \in \Omega \times \R$, all but a finite number of the series' terms vanish. This makes the series pointwise convergent and thus $Q$ measurable.
	\item For fixed $\theta \in \R$ and $A \subset \Omega$ open,
	\[
	\int_A P_\theta f (y) \d y = \int_{(T^\theta)^{-1} (A)} f (x) \d x 
	= \int_{(T^\theta)^{-1} (A)} \sum_{j \in \Z} \1_{U_j^\theta} (x) f(x) \d x
	\]
	whence we get by substituting $x = (T_{j}^{ \theta})^{-1} (y)$ in each term:
	\begin{align*}
		\int_A P_\theta f (y) \d y &= \int_A \sum_{j \in \Z} \left( \1_{U_j^\theta} \circ (T_{j}^{ \theta})^{-1} \right) \! (y) \left(f \circ (T_{j}^{ \theta})^{-1}\right) \! (y) \abs{\det D T_{j}^{\theta} \left((T_{j}^{ \theta})^{-1} (y)\right)} \d y \\
		&= \int_A \sum_{j \in \Z} \1_{T_{j}^{\theta} (U_j^\theta)} (y) \left(f \circ (T_{j}^{ \theta})^{-1}\right) \! (y) \abs{\frac{\partial \Phi_0}{\partial x_1} \left(\Gamma (T_{j}^{ \theta})^{-1} (y)\right)}^{-1} \d y \\
		&= \int_A Qf (y, \theta) \d y
	\end{align*}
	With this equality holding for open $A \subset \Omega$, a measure extension theorem allows us to induce that $P_\theta f(y) = Qf (y, \theta)$ for every $(y, \theta) \in \Omega \times \R$.
	\item $\bar Q$ is measurable as it is the product of two measurable functions, and for all $\theta \in \R$, $\abs{P_\theta f} \le P_\theta \abs{f}$ for which Fubini-Tonelli yields:
	\[
	\iint_{\Omega \times \R} H(\theta) P_\theta \abs{f}(y) \d\theta \d y = \int_\R H(\theta) \left( \int_\Omega P_\theta \abs{f} (y) \d y \right) \d\theta
	= \norm{f}_{L^1},
	\]
	hence the result.
	\qed
\end{enumerate}


This ensures that the operator described in \eqref{PF_moyen} is well-defined. Notice that the expression given for $P_\th$ in \eqref{forme_PF} allows us to establish an $L^\infty \to L^\infty$ bound:
\begin{equation}\label{borne_Linfini_PFmoy}
	\abs{P_\th f (x)} \le N \norm{f}_{L^\infty} C_1 ^{-\frac{k-1}{2}} C_2^{-\frac12} \sigma^{-\frac k2}
\end{equation}
where $N$ denotes the maximal number of nonempty open subsets $U_j^\th$ for fixed $\th \in \R$. Since this number has a bound that is uniform in $\th$ (see \eqref{nb-ouv}), we are able to derive a bound for the average operator on $L^\infty$ after establishing $L^1$-boundedness, which will be achieved by showing that $P$ is a Markov operator.

We recall that, for some measure space $(X, \T, \mu)$, a Markov operator over $X$ in the sense of Lasota \& Mackey (see \cite{LM}) is a linear operator $T : L^1 (X, \mu) \to L^1 (X, \mu)$ such that for all {\textit{nonnegative}} $f \in L^1$, $T f \ge 0$ and $\norm{T f}_{L^1} = \norm{f}_{L^1}$. 

\begin{theo}
	$P$ is a Markov operator.
\end{theo}
{\bf Proof} \\
Let $f \in L^1 (\Omega, m)$ such that $f \ge 0$. Thanks to  
properties (3.2.4), (3.2.5) from \cite{LM}, the $P_\th$ are Markov operators, hence $Pf = \int_\R H(\th) P_\th f \d\th \ge 0$, and 
\begin{align*}
	\norm{Pf}_{L^1} &= \int_\Omega Pf (y) \d y = \iint_{\Omega \times \R} H(\th) P_\th f (y) \d\th \d y \\
	&= \int_\R H(\th) \left(\int_\Omega P_\th f (y) \d y\right) \d\th = \int_\R H(\th)  \norm{P_\th f}_{L^1} \d\th \\
	&= \int_\R H(\th) \norm{f}_{L^1} \d\th = \norm{f}_{L^1},
\end{align*}
establishing that $P$ also is a Markov operator.
\qed

Proposition 3.1.1 of \cite{LM} then allows us to establish that:

\begin{cor}\label{borne_L1_PFmoy}
	For all $f \in L^1 (\Omega, m)$, $\norm{Pf}_{L^1} \le \norm{f}_{L^1}$.
\end{cor}

\begin{theo}
	For all $p \in [1, \infty]$, $P : L^p \to L^p$ is a bounded linear operator.
\end{theo}

{\textbf{Proof}} \\
By interpolation of \eqref{borne_Linfini_PFmoy} and Corollary \ref{borne_L1_PFmoy}, $P$ maps $L^p$ onto itself and by setting $s = 1 - \frac1p$, one has:
\[
\norm{P}_{\lin(L^p)} \le \norm{P}_{\lin(L^1)}^{1-s} \norm{P}_{\lin(L^\infty)}^s \le N^s C_1 ^{-\frac{(k-1)s}{2}} C_2^{-\frac s2} \sigma^{-\frac {ks}2},
\]
and thus $P \in \lin(L^p)$.
\qed

\subsection{Lasota-Yorke inequality and spectral properties}

Our aim now is to apply the Ionescu-Tulcea \& Marinescu spectral decomposition theorem on a suitable $P$-stable subspace. The space of functions of weak oscillation, as defined in \cite{ANV, SAU}, is one such subspace:
\begin{equation}
	V_1 (\Omega) = \left\{ f \in L^1 (\R^k) : \supp f \subset \Omega, \; \abs{f}_1 \coloneqq \sup_{0 < \eps < \eps_0} \eps^{-1} \int \osc\left(f, B_\eps(x)\right) \d x < +\infty \right\},
\end{equation}
where for any non-null measurable set $A$ and all $f \in L^1$, the oscillation of $f$ over $A$ is defined in \cite{ANV} by
\[
\osc(f, A) = \supess_{x, y \in A} \abs{f(x) - f(y)}.
\]
The first ingredient required for the decomposition theorem to apply is a Lasota-Yorke inequality on $P$. One such estimate was established in \cite{SAU} which can be applied to our $P_\th$:

\begin{lem}\label{D-eta-unif}
Under the conditions  \eqref{eta} to \eqref{cond} of Section \ref{sec-nch},
there exist $\eta\in [0,1[$ and $D\in \R^+$ such that, for all 
$\theta\in \R$, $\forall f\in V_1$, 
\begin{equation}\label{eq-D-eta-unif}
|P_{\theta}f|_1 \leq \eta |f|_1+ D \Vert f \Vert_{L^1_m}.
\end{equation}
\end{lem}
\noindent{\bf Proof.}\\
This is a consequence of the proof of Lemma 4.1 of \cite{SAU}. For 
$\eps \in (0,\eps_0) $ of Theorem \ref{unif},
 Lemma  4.1
gives a Lasota-Yorke inequality for the Perron-Frobenius operator $P_{\theta}$,
of the shape \eqref{eq-D-eta-unif}, with explicit $\eta$ and $D$ given by
$$
\eta= \left(1+ \frac{K}{\sqrt{\sigma}} \eps_0\right)
 \left(\frac{1}{\sqrt{\sigma}}+ 2 \sup\frac{G(\eps,\eps_0)\eps_0}{\eps}\right)
$$
and
$$
D = \frac{K}{\sqrt{\sigma}} \eps +2 \left(1+\frac{K}{\sqrt{\sigma}} \eps\right)
G(\eps,\eps_0).
$$
With the notation \eqref{etabarre}, we get 
$$
\eta \leq \left(1+ \frac{K}{\sqrt{\sigma}} \eps_0\right)\overline{\eta(\eps_0)},\quad
D \leq \frac{K}{\sqrt{\sigma}} \eps_0 + 
\left(1+\frac{K}{\sqrt{\sigma}} \eps_0\right)\overline{\eta(\eps_0)}.
$$
For $\eps_0$ sufficiently small, $\eta$ is smaller than $1$. 
Neither of them depends on $\theta$ since the constants $K,\sigma$ and the 
function $\overline{\eta(\eps_0)}$ do not depend on $\theta$.
\qed 

It follows from Lemma \ref{D-eta-unif} that $V_1$ is invariant under the operators $P_\th$ for all $\th \in \R$. In the following, we shall establish a similar estimate for $P$.

\begin{lem}\label{measur_osc}
	For all $f \in V_1 (\Omega)$ and $\eps > 0$, the map 
	\[
	\left.\begin{matrix}
		\R^k \times \R & \longrightarrow & \R \\
		(x, \theta) & \longmapsto & \osc(P_\theta f, B_\eps (x))
	\end{matrix}\right.
	\]
	is measurable.
\end{lem}

{\bf Proof} \\
Let $a > 0$ and $O_\th : (y, z) \in \Omega^2 \mapsto \abs{P_\th f (y) - P_\th f (z)}$. We shall denote by $\chi_a$ the indicator function of the set
\[
\left\{ (x, y, z, \th) \in \R^k \times \Omega \times \Omega \times \R : y, z \in B_\eps (x) \text{ and } O_\th (y, z) > a \right\},
\]
whose measurability is easily deduced from that of $(\th, y)
\mapsto P_\th f(y)$.
Define
\[
M_a : (x, \th) \in \R^k \times \R \mapsto \iint_{\Omega^2} \chi_a (x, y, z, \th) \d y \d z
\]
By the Fubini-Tonelli theorem, $M_a$ is measurable, thus so is the set
\begin{align*}
	\left\{(x, \th) \in \R^k \times \R : \supess_{y, z \in B_\eps (x)} O_\th (y, z) > a\right\} &= \left\{(x, \th) : m\left(\left\{(y, z) \in B_\eps (x)^2 : O_\th(y, z) > a\right\}\right) > 0\right\} \\
	&= \left\{(x, \th) : M_a (x, \th) > 0\right\} = M_a^{-1} \left( (0, +\infty) \right)
\end{align*}
Hence the map $\displaystyle (x, \th) \in \R^k \times \R \mapsto \supess_{y, z \in B_\eps (x)} O_\th (y, z)$ is measurable, which concludes the proof.
\qed

The proof for the main result in this section utilizes the following well-known estimate: 
\[
\sup_{x \in X} \int_Y g(x, y) \d\mu_Y (y) \le \int_Y \sup_{x \in X} g(x, y) \d\mu_Y (y)
\]
for all nonnegative, bounded $g \in L^1 (X \times Y, \mu_X \otimes \mu_Y)$. A measure-theoretic version shall also be utilized, for which we give a short proof.

\begin{lem}
	Let $(X, \T_X, \mu_X)$, $(Y, \T_Y, \mu_Y)$ be $\sigma$-finite measure spaces, and consider a nonnegative essentially bounded $g \in L^1 (X \times Y, \mu_X \otimes \mu_Y)$. The mapping $y \in Y \mapsto \supess_{x \in A} g(x, y)$ is measurable, and
	\begin{equation}\label{ineg_supess_int}
		\supess_{x \in A} \int_Y g(x, y) \d\mu_Y (y) \le \int_Y \supess_{x \in A} g(x, y) \d\mu_Y (y)
	\end{equation}
	for all $A \in \T_X$ with positive measure.
\end{lem}

{\bf Proof} \\
Measurability is established as in
 Lemma \ref{measur_osc}. Consider $f : (x, y) \mapsto \supess_A g(\cdot, y) - g(x, y)$ and denote by $\chi$ the indicator of $f^{-1} (\R^*_-)$. By Fubini-Tonelli,
\[
\iint_{X \times Y} \chi = \int_Y \mu_X \left( \left\{ x \in X : g(x, y) > \supess_A g(\cdot, y) \right\} \right) \d \mu_Y (y) = 0,
\]
which entails that $f^{-1} (\R^*_-)$ is null and thus that $f \ge 0$ almost everywhere. Therefore, for $\mu_X$-almost every $x \in X$,
\[
\int_Y g(x, y) \d\mu_Y (y) \le \int_Y \supess_A g(\cdot, y) \d\mu_Y (y),
\]
hence \eqref{ineg_supess_int}.
\qed

We may now move on to the proof of this section's main result:

\begin{prop}
There exist $\eta \in (0, 1)$ and $C > 0$ such that,
for all $f \in V_1 (\Omega)$, $Pf\in V_1$ and 
	\[
	\norm{Pf}_{V_1} \le \eta \norm{f}_{V_1} + C \norm{f}_{L^1}
	\]
\end{prop}
{\bf Proof} \\
By Lemma \ref{D-eta-unif}, for all $f \in V_1$ and $\th \in \R$, 
\[
\abs{P_\th f}_1 \le \eta \abs{f}_1 + D \norm{f}_{L^1}
\]
where $\eta\in \left(0, 1\right)$, $D > 0$ do not depend on $f$ or $\th$. Furthermore, \eqref{ineg_supess_int} and the definition of $P$ yield:
{
	\allowdisplaybreaks
	\begin{align*}
		\abs{Pf}_1 &= \sup_{0 < \eps < \eps_0} \eps^{-1} \int_{\R^k} \supess_{y, z \in B_\eps (x)} \abs{Pf(y) - Pf(z)} \d x \\
		&= \sup_{0 < \eps < \eps_0} \eps^{-1} \int_{\R^k} \supess_{y, z \in B_\eps (x)} \abs{\int_\R H(\th) \left(P_\th f(y) - P_\th f(z)\right) \d\th} \d x \\
		&\le \sup_{0 < \eps < \eps_0} \eps^{-1} \int_{\R^k} \supess_{y, z \in B_\eps (x)} \int_\R H(\th) \abs{P_\th f(y) - P_\th f(z)} \d\th \d x \\
		&\le \sup_{0 < \eps < \eps_0} \eps^{-1} \int_{\R^k} \int_\R H(\th) \supess_{y, z \in B_\eps (x)} \abs{P_\th f(y) - P_\th f(z)} \d\th \d x \\
		&= \sup_{0 < \eps < \eps_0} \eps^{-1} \int_\R \int_{\R^k} H(\th) \osc \left(P_\th f, B_\eps (x) \right) \d x \d \th \\
		&\le \int_\R H(\th) \left( \sup_{0 < \eps < \eps_0} \eps^{-1} \int_{\R^k} \osc \left(P_\th f, B_\eps (x)\right) \d x \right) \d \th = \int_\R H(\th) \abs{P_\th f}_1 \d \th,
	\end{align*}
}
and it follows from \eqref{eq-D-eta-unif} that
\begin{equation}\label{las-yor-prelim}
	\abs{Pf}_1 \le \int_\R H(\th) \left(\eta \abs{f}_1 + D \norm{f}_{L^1}\right) \d \th = \eta \abs{f}_1 + D \norm{f}_{L^1} .
\end{equation}
Recall in addition that $\norm{Pf}_{L^1} \le \norm{f}_{L^1}$ for all $f \in V_1$ ; adding this to \eqref{las-yor-prelim} yields
\begin{align*}
	\norm{Pf}_{V_1} &\le \eta \abs{f}_1 + D \norm{f}_{L^1} + \norm{f}_{L^1} \\
	&= \eta \abs{f}_1 + D \norm{f}_{L^1} + (\eta + 1 - \eta) \norm{f}_{L^1} \\
	&= \eta \norm{f}_{V_1} + (D + 1 - \eta) \norm{f}_{L^1},
\end{align*}
and the result is established for $C \coloneqq D + 1 - \eta$.
\qed

%

We now have every ingredient we need for Ionescu-Tulcea \& Marinescu's
 theorem to apply, which proves points \ref{projs}, 
\ref{spec-gap}, and the beginning of point \ref{eigs} of Theorem \ref{av-TPF}. 
The operator $P$ decomposes as 
\[
P = S + \sum_{j = 1}^p \lambda_j \Pi_j,
\]
with the commutation relations of Theorem \ref{av-TPF}. The fact that $\lambda_1 = 1$ is proved in the following Lemma.

\begin{lem}
	The averaged transfer operator $P$ has the eigenvalue $\lambda_1 = 1$, and $h_* \coloneqq \Pi_1 \left(\frac1{m(\Omega)}\right)$ is a nonnegative fixed point.
\end{lem}

{\bf Proof}\\
Applying Lemma (4.1) of \cite{ITM} with $c=1$, for all $f\in V_1$
gives a function $\tilde{f}\in V_1$ such that 
$$
\lim_{n\rightarrow \infty} 
\norm{\tilde{f} - \frac{1}{n}\sum_{k=1}^n P^k f}_{L^1} =0. 
$$
Taking $\lambda_1 = 1$, the projector $\Pi_1$ from Theorem \ref{av-TPF} on $V_1$ is then defined 
by setting $\tilde{f}= \Pi_1 f$. 
It is a positive operator, in the sense that, if $g\in V_1$ and $g\geq 0$, then
$\Pi_1 g\geq 0$ too (consider an a.e. converging subsequence of
$ \frac{1}{n}\sum_{k=1}^n  P^k g $,
which is a sequence of $\geq 0$ functions). 

One may check that 
$$
\int_{\Omega} \frac{1}{n} \sum_{k=1}^n  P^k f \, \d m  = \int_{\Omega} f \, \d m,
$$
since $\int_{\Omega}P f \, \d m = \int_{\Omega} f \, \d m$ due to $P$ being Markov. 
Hence, if $f>0$, $\Pi_1 f$ is not the zero function because
$  \int_{\Omega} \Pi_1 f \, \d m = \int_{\Omega} f \, \d m\neq 0$. 
Moreover, the function $\tilde f= \Pi_1 f$  satisfies 
$P \tilde f = \tilde f$.
\\
Thus $h_*= \Pi_1(\frac1{m(\Omega)})$ is a fixed point of $P$ and
it is nonnegative by positiveness of $\Pi_1$.
\qed

{\bf Proof of Theorem \ref{th-inv-meas}}

Here $\theta= (\theta_n)_{n\in \N}=(\theta_0,\check{\theta}) \in \R^{\N}$ is a
real sequence. In \eqref{Effe}, we defined a mapping $F$ on
$\R^{\N}\times \Omega$,
 setting $F((\theta_0,\check{\theta}),x)= (\check{\theta},T^{\theta_0}(x))$.
We seek to establish that $\Pr_\Theta \otimes \mu$ is invariant for $F$.  It suffices to show that
\[
(\Pr_\Theta \otimes \mu)(A \times B) = (\Pr_\Theta \otimes \mu) \left( F^{-1} (A \times B) \right)
\]
for all $A = A_0 \times ... \times A_n \times \R \times ...$ with $A_0, ..., A_n \in \bor(\R)$, $n \in \N$, and all $B \in \bor(\Omega)$, as sets of this form generate the Borel $\sigma$-algebra over $\R^\N \times \Omega$.

We have
\begin{align*}
	(\Pr_\Theta \otimes \mu) \left( F^{-1} (A \times B) \right) &= \iint_{\R^\N \times \Omega} \1_{A \times B} \left(F(\theta, x)\right) \d (\Pr_\Theta \otimes \mu) (\theta, x) \\
	&= \int_{\R^\N} \1_A (S\theta) \left[ \int_\Omega \1_B \left( T^{\theta_0} (x) \right) \d\mu(x) \right] \d \Pr_\Theta (\theta) \\
	&= \iint_{\R\times\R^\N} \1_{\R \times A} (\theta_0, \check \theta) \left[ \int_{(T^{\theta_0})^{-1} (B)} \d\mu(x) \right] \d (\Pr_{\Theta_0} \otimes \Pr_\Theta) (\theta_0, \check \theta),
	\intertext{thanks to the independence of the variables $(\Theta_t)$, which allows us to split the integral as follows:}
	(\Pr_\Theta \otimes \mu) \left( F^{-1} (A \times B) \right) &= \int_A \left[ \int_\R \int_{(T^{\theta_0})^{-1} (B)} \d\mu(x) \d\Pr_{\Theta_0} (\theta_0) \right] \d \Pr_\Theta (\check \theta) \\
	&= \int_A \left[ \int_\R \int_{(T^{\theta_0})^{-1} (B)} h_* \d m(x) \d\Pr_{\Theta_0} (\theta_0) \right] \d \Pr_\Theta (\check \theta) \\
	&= \int_A \left[ \int_\R \int_B P_{\theta_0} h_* \d m(x) \d\Pr_{\Theta_0} (\theta_0) \right] \d \Pr_\Theta (\check \theta)
	\intertext{by definition of $\mu$ and $P_{\theta_0}$. It then follows by applying Fubini's theorem that}
	(\Pr_\Theta \otimes \mu) \left( F^{-1} (A \times B) \right) &= \int_A \left[ \int_B \int_\R P_{\theta_0} h_* \d\Pr_{\Theta_0} (\theta_0) \d m(x) \right] \d \Pr_\Theta (\check \theta) \\
	&= \int_A \left[ \int_B P h_* \d m(x) \right] \d \Pr_\Theta (\check \theta),
	\intertext{whence, since $Ph_* = h_*$,}
	(\Pr_\Theta \otimes \mu) \left( F^{-1} (A \times B) \right) &= \int_A \left[ \int_B \d \mu(x) \right] \d \Pr_\Theta (\check \theta) \\
	&= \iint_{A \times B} \d (\Pr_\Theta \otimes \mu) (\check \theta, x) = (\Pr_\Theta \otimes \mu) (A \times B),
\end{align*}
which establishes the invariance of $\mu$. 
\qed

\section*{Declarations}
The authors have no competing interests to declare.

\end{document}